\documentclass[10pt]{amsart}
\usepackage{amsmath}
\usepackage{amsxtra}
\usepackage{amscd}
\usepackage{amsthm}
\usepackage{amsfonts}
\usepackage{amssymb}
\usepackage{eucal}
\usepackage{scalerel}
\usepackage[matrix,arrow,curve]{xy}

\usepackage[
pdftex,
bookmarks=false,
colorlinks=true,
debug=true,
pdfnewwindow=true]{hyperref}

\theoremstyle{definition}

\theoremstyle{remark}

\newcommand{\nc}{\newcommand}
\nc{\renc}{\renewcommand} \nc{\ssec}{\subsection}
\nc{\sssec}{\subsubsection} 

\nc{\on}{\operatorname} \nc{\wh}{\widehat}
\nc\ol{\overline} \nc\ul{\underline} \nc\wt{\widetilde}

\nc{\BA}{{\mathbb{A}}} \nc{\BC}{{\mathbb{C}}} \nc{\BQ}{{\mathbb{Q}}}
\nc{\BM}{{\mathbb{M}}} \nc{\BN}{{\mathbb{N}}}
\nc{\BP}{{\mathbb{P}}} \nc{\BR}{{\mathbb{R}}}
\nc{\BZ}{{\mathbb{Z}}} \nc{\BS}{{\mathbb{S}}} \nc{\BW}{{\mathbb{W}}}

\nc{\CA}{{\mathcal{A}}} \nc{\CB}{{\mathcal{B}}}
\nc{\CE}{{\mathcal{E}}} \nc{\CF}{{\mathcal{F}}}
\nc{\CG}{{\mathcal{G}}} \nc{\CH}{{\mathcal{H}}}
\nc{\CI}{{\mathcal{I}}} \nc{\CK}{{\mathcal{K}}} \nc{\CL}{{\mathcal{L}}}
\nc{\CM}{{\mathcal{M}}} \nc{\CN}{{\mathcal{N}}}
\nc{\CO}{{\mathcal{O}}} \nc{\CP}{{\mathcal{P}}}
\nc{\CQ}{{\mathcal{Q}}} \nc{\CR}{{\mathcal{R}}}
\nc{\CS}{{\mathcal{S}}} \nc{\CT}{{\mathcal{T}}}
\nc{\CU}{{\mathcal{U}}} \nc{\CV}{{\mathcal{V}}}  \nc{\CY}{{\mathcal Y}}
\nc{\CW}{{\mathcal{W}}} \nc{\CZ}{{\mathcal{Z}}}

\nc{\cM}{{\check{\mathcal M}}{}} \nc{\csM}{{\check{\mathcal A}}{}}
\nc{\oM}{{\overset{\circ}{\mathcal M}}{}}
\nc{\obM}{{\overset{\circ}{\mathbf M}}{}}
\nc{\oCA}{{\overset{\circ}{\mathcal A}}{}}
\nc{\obA}{{\overset{\circ}{\mathbf A}}{}}
\nc{\ooM}{{\overset{\circ}{M}}{}}
\nc{\osM}{{\overset{\circ}{\mathsf M}}{}}
\nc{\vM}{{\overset{\bullet}{\mathcal M}}{}}
\nc{\nM}{{\underset{\bullet}{\mathcal M}}{}}
\nc{\oD}{{\overset{\circ}{\mathcal D}}{}}
\nc{\obD}{{\overset{\circ}{\mathbf D}}{}}
\nc{\oA}{{\overset{\circ}{\mathbb A}}{}}
\nc{\op}{{\overset{\bullet}{\mathbf p}}{}}
\nc{\cp}{{\overset{\circ}{\mathbf p}}{}}
\nc{\oU}{{\overset{\bullet}{\mathcal U}}{}}
\nc{\ofZ}{{\overset{\circ}{\mathfrak Z}}{}}

\nc{\ff}{{\mathfrak{f}}} \nc{\fv}{{\mathfrak{v}}}
\nc{\fa}{{\mathfrak{a}}} \nc{\fb}{{\mathfrak{b}}}
\nc{\fd}{{\mathfrak{d}}} \nc{\fe}{{\mathfrak{e}}}
\nc{\fg}{{\mathfrak{g}}} \nc{\fgl}{{\mathfrak{gl}}}
\nc{\fh}{{\mathfrak{h}}} \nc{\fri}{{\mathfrak{i}}}
\nc{\fj}{{\mathfrak{j}}} \nc{\fk}{{\mathfrak{k}}}
\nc{\fm}{{\mathfrak{m}}} \nc{\fn}{{\mathfrak{n}}}
\nc{\ft}{{\mathfrak{t}}} \nc{\fu}{{\mathfrak{u}}}
\nc{\fw}{{\mathfrak{w}}} \nc{\fz}{{\mathfrak{z}}}
\nc{\fp}{{\mathfrak{p}}} \nc{\frr}{{\mathfrak{r}}}
\nc{\fs}{{\mathfrak{s}}} \nc{\fsl}{{\mathfrak{sl}}}
\nc{\hsl}{{\widehat{\mathfrak{sl}}}}
\nc{\hgl}{{\widehat{\mathfrak{gl}}}}
\nc{\hg}{{\widehat{\mathfrak{g}}}}
\nc{\chg}{{\widehat{\mathfrak{g}}}{}^\vee}
\nc{\hn}{{\widehat{\mathfrak{n}}}}
\nc{\chn}{{\widehat{\mathfrak{n}}}{}^\vee}

\nc{\fA}{{\mathfrak{A}}} \nc{\fB}{{\mathfrak{B}}}
\nc{\fD}{{\mathfrak{D}}} \nc{\fE}{{\mathfrak{E}}}
\nc{\fF}{{\mathfrak{F}}} \nc{\fG}{{\mathfrak{G}}} \nc{\fH}{{\mathfrak{H}}}
\nc{\fI}{{\mathfrak{I}}} \nc{\fJ}{{\mathfrak{J}}}
\nc{\fK}{{\mathfrak{K}}} \nc{\fL}{{\mathfrak{L}}}
\nc{\fM}{{\mathfrak{M}}} \nc{\fN}{{\mathfrak{N}}}
\nc{\frP}{{\mathfrak{P}}} \nc{\fQ}{{\mathfrak{Q}}}
\nc{\fT}{{\mathfrak{T}}} \nc{\fU}{{\mathfrak{U}}}
\nc{\fV}{{\mathfrak{V}}} \nc{\fW}{{\mathfrak{W}}}
\nc{\fX}{{\mathfrak{X}}} \nc{\fY}{{\mathfrak{Y}}}
\nc{\fZ}{{\mathfrak{Z}}}

\nc{\ba}{{\mathbf{a}}}
\nc{\bb}{{\mathbf{b}}} \nc{\bc}{{\mathbf{c}}}
\nc{\be}{{\mathbf{e}}} \nc{\bj}{{\mathbf{j}}}
\nc{\bn}{{\mathbf{n}}} \nc{\bp}{{\mathbf{p}}}
\nc{\bq}{{\mathbf{q}}} \nc{\br}{{\mathbf{r}}} \nc{\bt}{{\mathbf{t}}}
\nc{\bfu}{{\mathbf{u}}} \nc{\bv}{{\mathbf{v}}}
\nc{\bx}{{\mathbf{x}}} \nc{\by}{{\mathbf{y}}}
\nc{\bw}{{\mathbf{w}}} \nc{\bA}{{\mathbf{A}}}
\nc{\bB}{{\mathbf{B}}} \nc{\bC}{{\mathbf{C}}}
\nc{\bD}{{\mathbf{D}}} \nc{\bF}{{\mathbf{F}}}
\nc{\bH}{{\mathbf{H}}} \nc{\bJ}{{\mathbf{J}}} \nc{\bK}{{\mathbf{K}}}
\nc{\bM}{{\mathbf{M}}} \nc{\bN}{{\mathbf{N}}}
\nc{\bO}{{\mathbf{O}}} \nc{\bS}{{\mathbf{S}}} \nc{\bT}{{\mathbf{T}}}
\nc{\bV}{{\mathbf{V}}} \nc{\bW}{{\mathbf{W}}}
\nc{\bX}{{\mathbf{X}}}
\nc{\bY}{{\mathbf{Y}}} \nc{\bP}{{\mathbf{P}}}
\nc{\bZ}{{\mathbf{Z}}} \nc{\bh}{{\mathbf{h}}}

\nc{\sA}{{\mathsf{A}}} \nc{\sB}{{\mathsf{B}}}
\nc{\sC}{{\mathsf{C}}} \nc{\sD}{{\mathsf{D}}}
\nc{\sE}{{\mathsf{E}}} \nc{\sF}{{\mathsf{F}}} \nc{\sG}{{\mathsf{G}}}
\nc{\sI}{{\mathsf{I}}} \nc{\sK}{{\mathsf{K}}} \nc{\sL}{{\mathsf{L}}}
\nc{\sM}{{\mathsf{M}}} \nc{\sO}{{\mathsf{O}}}
\nc{\sQ}{{\mathsf{Q}}} \nc{\sP}{{\mathsf{P}}}
\nc{\sT}{{\mathsf{T}}} \nc{\sZ}{{\mathsf{Z}}}
\nc{\sV}{{\mathsf{V}}} \nc{\sW}{{\mathsf{W}}}
\nc{\sfp}{{\mathsf{p}}} \nc{\sr}{{\mathsf{r}}}
\nc{\st}{{\mathsf{t}}} \nc{\sfb}{{\mathsf{b}}}
\nc{\sfc}{{\mathsf{c}}} \nc{\sd}{{\mathsf{d}}}
\nc{\sz}{{\mathsf{z}}}

\nc{\tA}{{\widetilde{\mathbf{A}}}}
\nc{\tB}{{\widetilde{\mathcal{B}}}}
\nc{\tg}{{\widetilde{\mathfrak{g}}}} \nc{\tG}{{\widetilde{G}}}
\nc{\TM}{{\widetilde{\mathbb{M}}}{}}
\nc{\tO}{{\widetilde{\mathsf{O}}}{}}
\nc{\tU}{{\widetilde{\mathfrak{U}}}{}} \nc{\TZ}{{\tilde{Z}}}
\nc{\tx}{{\tilde{x}}} \nc{\tbv}{{\tilde{\bv}}}
\nc{\tfP}{{\widetilde{\mathfrak{P}}}{}} \nc{\tz}{{\tilde{\zeta}}}
\nc{\tmu}{{\tilde{\mu}}}

\nc{\urho}{\underline{\rho}} \nc{\uB}{\underline{B}}
\nc{\uC}{{\underline{\mathbb{C}}}} \nc{\ui}{\underline{i}}
\nc{\uj}{\underline{j}} \nc{\ofP}{{\overline{\mathfrak{P}}}}
\nc{\oB}{{\overline{\mathcal{B}}}}
\nc{\og}{{\overline{\mathfrak{g}}}} \nc{\oI}{{\overline{I}}}

\nc{\eps}{\varepsilon} \nc{\hrho}{{\hat{\rho}}}
\nc{\blambda}{{\bar\lambda}} \nc{\bmu}{{\bar\mu}} \nc{\bnu}{{\bar\nu}}

\nc{\one}{{\mathbf{1}}} \nc{\two}{{\mathbf{t}}}

\nc{\Sym}{{\mathop{\operatorname{\rm Sym}}}}
\nc{\Tot}{{\mathop{\operatorname{\rm Tot}}}}
\nc{\Spec}{{\mathop{\operatorname{\rm Spec}}}}
\nc{\Ker}{{\mathop{\operatorname{\rm Ker}}}}
\nc{\Hilb}{{\mathop{\operatorname{\rm Hilb}}}}
\nc{\End}{{\mathop{\operatorname{\rm End}}}}
\nc{\Ext}{{\mathop{\operatorname{\rm Ext}}}}
\nc{\Hom}{{\mathop{\operatorname{\rm Hom}}}}
\nc{\CHom}{{\mathop{\operatorname{{\mathcal{H}}\it om}}}}
\nc{\GL}{{\mathop{\operatorname{\rm GL}}}}
\nc{\gr}{{\mathop{\operatorname{\rm gr}}}}
\nc{\Id}{{\mathop{\operatorname{\rm Id}}}}
\nc{\defi}{{\mathop{\operatorname{\rm def}}}}
\nc{\length}{{\mathop{\operatorname{\rm length}}}}
\nc{\supp}{{\mathop{\operatorname{\rm supp}}}}

\nc{\Cliff}{{\mathsf{Cliff}}}
\nc{\Fl}{{\mathsf{Fl}}} \nc{\Fib}{{\mathsf{Fib}}}
\nc{\Coh}{{\mathsf{Coh}}} \nc{\FCoh}{{\mathsf{FCoh}}}

\nc{\reg}{{\text{\rm reg}}}

\nc{\cplus}{{\mathbf{C}_+}} \nc{\cminus}{{\mathbf{C}_-}}
\nc{\cthree}{{\mathbf{C}_*}} \nc{\Qbar}{{\bar{Q}}}

\newcommand{\oZ}{\vphantom{j^{X^2}}\smash{\overset{\circ}{\vphantom{\rule{0pt}{0.55em}}\smash{Z}}}}
\newcommand\iso{\,\vphantom{j^{X^2}}\smash{\overset{\sim}{\vphantom{\rule{0pt}{0.20em}}\smash{\longrightarrow}}}\,}
\nc{\Gtimes}{\vphantom{j^{X^2}}\smash{\overset{\mathsf G}{\vphantom{\rule{0pt}{0.30em}}\smash{\times}}}}

\nc{\bOmega}{{\overline{\Omega}}}

\nc{\seq}[1]{\stackrel{#1}{\sim}}

\nc{\aff}{{\operatorname{aff}}}
\nc{\fin}{{\operatorname{fin}}}
\nc{\Gr}{{\operatorname{Gr}}}
\nc{\GR}{{\mathbf{Gr}}}
\nc{\Perv}{{\operatorname{Perv}}}
\nc{\Rep}{{\operatorname{Rep}}}
\nc{\IC}{{\operatorname{IC}}}
\nc{\Bun}{{\operatorname{Bun}}}
\nc{\Proj}{{\operatorname{Proj}}}
\nc{\pt}{{\operatorname{pt}}}
\nc{\bfmu}{{\boldsymbol{\mu}}}

\begin{document}

\author{Michael Finkelberg}
\title
[Double affine Grassmannians and Coulomb branches of quiver gauge theories]
{Double affine Grassmannians and Coulomb branches of $3d\ \CN=4$ quiver gauge
theories}





\address{{\it Address}:\newline
National Research University Higher School of Economics, Russian Federation\newline
Department of Mathematics, 6 Usacheva st., Moscow 119048;\newline
Skolkovo Institute of Science and Technology;\newline
Institute for Information Transmission Problems}

\email{\newline fnklberg@gmail.com}

\begin{abstract}
We propose a conjectural construction of various slices for double affine
Grassmannians as Coulomb branches of 3-dimensional $\CN=4$ supersymmetric
affine quiver gauge theories. It generalizes the known construction for the
usual affine Grassmannians, and makes sense for arbitrary symmetric Kac-Mody
algebras.
\end{abstract}
\maketitle

\section{Introduction}

\subsection{Historical background}
The geometric Satake equivalence~\cite{lu,gi,bd,mv} proposed by V.~Drinfeld 
for the needs of
the Geometric Langlands Program proved very useful for the study of 
representation theory of reductive algebraic groups (starting from G.~Lusztig's
construction of $q$-analogues of weight multiplicities). About 15 years ago,
I.~Frenkel and I.~Grojnowski envisioned an extension of the geometric Satake
equivalence to the case of loop groups. The affine Grassmannians (the main
objects of the geometric Satake equivalence) are ind-schemes of ind-finite type.
Their loop analogues (double affine Grassmannians) are much more infinite,
beyond our current technical abilities. We are bound to settle for some
provisional substitutes, such as transversal slices to the smaller strata
in the closures of bigger strata. These substitutes still carry quite  
powerful geometric information.

Following I.~Frenkel's suggestion, some particular slices for the double
affine Grassmannians were constructed in terms of Uhlenbeck compactifications
of instanton moduli spaces on Kleinian singularities about 10 years ago.
More recently, H.~Nakajima's approach to Coulomb branches of 3-dimensional
$\CN=4$ supersymmetric gauge theories, applied to affine quiver gauge theories,
paved a way for the construction of the most general slices.

\subsection{Contents}
We recall the geometric Satake equivalence in~Section~\ref{satake}.
The (generalized) slices for the affine Grassmannians are reviewed
in~Section~\ref{slices}. The problem of constructing (the slices for) the
double affine Grassmannians is formulated in~Section~\ref{dag}. 
The mathematical construction of Coulomb branches of $3d\ \CN=4$ gauge theories
and its application to slices occupies~Section~\ref{coulomb}.
Some more applications are mentioned in~Section~\ref{applications}.

\subsection{Acknowledgments}
This report is based mostly on the works of A.~Braverman and H.~Nakajima, some 
joint with the author. I was incredibly lucky to have an opportunity to learn 
mathematics from them. Before meeting them, I was introduced to some 
semiinfinite ideas sketched below by A.~Beilinson, V.~Drinfeld, B.~Feigin,
V.~Ginzburg and I.~Mirkovi\'c. It is also a pleasure to acknowledge my 
intellectual debt to R.~Bezrukavnikov, D.~Gaiotto, D.~Gaitsgory, J.~Kamnitzer 
and V.~Pestun.

\section{Geometric Satake equivalence}
\label{satake}
Let $\CO$ denote the formal power series ring $\BC[[z]]$, and let $\CK$ denote
its fraction field $\BC((z))$. Let $G$ be an almost simple complex algebraic
group with a Borel and a Cartan subgroup $G\supset B\supset T$, and with
the Weyl group $W_\fin$ of $(G,T)$. Let $\Lambda$ 
be the coweight lattice, and let $\Lambda^+\subset\Lambda$ be the submonoid of
dominant coweights. Let also $\Lambda_+\subset\Lambda$ be the submonoid spanned
by the simple coroots $\alpha_i,\ i\in I$.
We denote by $G^\vee\supset T^\vee$ the Langlands dual
group, so that $\Lambda$ is the weight lattice of $G^\vee$.

The affine Grassmannian $\Gr_G=G_\CK/G_\CO$ is an ind-projective scheme,
the union $\bigsqcup_{\blambda\in\Lambda^+}\Gr_G^\blambda$ of $G_\CO$-orbits.
The closure of $\Gr_G^\blambda$ is a projective variety 
$\ol\Gr{}^\blambda_G=\bigsqcup_{\bmu\leq\blambda}\Gr^\bmu_G$. The fixed point set
$\Gr^T_G$ is naturally identified with the coweight lattice $\Lambda$;
and $\bmu\in\Lambda$ lies in $\Gr_G^\lambda$ iff $\bmu\in W_\fin\blambda$.

One of the cornerstones of the Geometric Langlands Program initiated by 
V.~Drinfeld is an equivalence $\BS$ of the tensor category $\Rep(G^\vee)$ and
the category $\Perv_{G_\CO}(\Gr_G)$ of $G_\CO$-equivariant perverse constructible
sheaves on $\Gr_G$ equipped with a natural monoidal convolution structure 
$\star$ and a fiber functor $H^\bullet(\Gr_G,-)$~\cite{lu,gi,bd,mv}. It is a 
categorification of the classical Satake isomorphism between 
$K(\Rep(G^\vee))=\BC[T^\vee]^{W_\fin}$ and the spherical affine Hecke algebra of 
$G$. The geometric Satake equivalence $\BS$ sends an irreducible $G^\vee$-module
$V^\blambda$ with highest weight $\blambda$ to the Goresky-MacPherson sheaf
$\IC(\ol\Gr{}^\blambda_G)$.

In order to construct a commutativity constraint for 
$(\Perv_{G_\CO}(\Gr_G),\star)$, Beilinson and Drinfeld introduced a relative
version $\Gr_{G,BD}$ of the Grassmannian over the Ran space of a smooth curve 
$X$, and a fusion monoidal structure $\Psi$ on $\Perv_{G_\CO}(\Gr_G)$ (isomorphic 
to $\star$). One of the main discoveries of~\cite{mv} was a $\Lambda$-grading
of the fiber functor 
$H^\bullet(\Gr_G,\CF)=\bigoplus_{\blambda\in\Lambda}\Phi_\blambda(\CF)$ by the 
hyperbolic stalks at $T$-fixed points. For a $G^\vee$-module $V$, its weight
space $V_\blambda$ is canonically isomorphic to the hyperbolic stalk
$\Phi_\blambda(\BS V)$.

Various geometric structures of a perverse sheaf $\BS V$ reflect some fine 
representation theoretic structures of $V$, such as Brylinski-Kostant filtration
and the action of dynamical Weyl group, see~\cite{gr}. One of the important 
technical tools of studying $\Perv_{G_\CO}(\Gr_G)$ is the embedding 
$\Gr_G\hookrightarrow\GR_G$ into Kashiwara infinite type scheme
$\GR_G=G_{\BC((z^{-1}))}/G_{\BC[z]}$~\cite{ka1,kt}. The quotient 
$G_{\BC[[z^{-1}]]}\backslash\GR_G$ is the moduli stack $\Bun_G(\BP^1)$ of 
$G$-bundles on the projective line $\BP^1$.
The $G_{\BC[[z^{-1}]]}$-orbits on $\GR_G$ are of finite codimension; they are also
numbered by the dominant coweights of $G$, and the image of an orbit 
$\GR_G^\blambda$ in $\Bun_G(\BP^1)$ consists of $G$-bundles of isomorphism type
$\blambda$~\cite{gro}. The stratifications 
$\Gr_G=\bigsqcup_{\blambda\in\Lambda^+}\Gr_G^\blambda$ and
$\GR_G=\bigsqcup_{\blambda\in\Lambda^+}\GR_G^\blambda$ are transversal, and their
intersections and various generalizations thereof are the subject of the next 
section.

\section{Generalized slices}
\label{slices}

\subsection{The dominant case}
\label{dominant}
We denote by $K_1$ the first congruence subgroup of $G_{\BC[[z^{-1}]]}$: the kernel
of the evaluation projection 
$\on{ev}_\infty\colon G_{\BC[[z^{-1}]]}\twoheadrightarrow G$. The transversal
slice $\CW_\bmu^\blambda$ (resp.\ $\ol\CW{}_\bmu^\blambda$) is defined as the 
intersection of $\Gr_G^\blambda$ (resp.\ $\ol\Gr{}_G^\blambda$) and 
$K_1\cdot\bmu$ in $\GR_G$. It is known that $\ol\CW{}_\bmu^\blambda$ is nonempty
iff $\bmu\leq\blambda$, and $\dim\ol\CW{}_\bmu^\blambda$ is an affine irreducible
variety of dimension $\langle2\bar\rho^{\!\scriptscriptstyle\vee},\blambda-\bmu\rangle$.
Following an idea of I.~Mirkovi\'c,~\cite{kwy} proved that 
$\ol\CW{}_\bmu^\blambda=\bigsqcup_{\bmu\leq\bar{\nu}\leq\blambda}\CW_\bmu^{\bar\nu}$
is the decomposition of $\ol\CW{}_\bmu^\blambda$ into symplectic leaves of a
natural Poisson structure. 

The only $T$-fixed point of $\ol\CW{}_\bmu^\blambda$ is $\bmu$. We consider the
cocharacter $2\bar\rho\colon\BC^\times\to T$, and denote by 
$R^\blambda_\bmu\subset\ol\CW{}_\bmu^\blambda$ the corresponding repellent: the
closed affine subvariety formed by all the points that flow into $\bmu$
under the action of $2\bar\rho(t)$, as $t$ goes to $\infty$.
Let $r$ stand for the closed embedding of $R^\blambda_\bmu$ into 
$\ol\CW{}_\bmu^\blambda$, and let $\iota$ stand for the closed embedding of
$\bmu$ into $R^\blambda_\bmu$. Then the hyperbolic stalk $\Phi^\blambda_\bmu\CF$
of a $T$-equivariant constructible complex $\CF$ on $\ol\CW{}_\bmu^\blambda$
is defined as $\iota^!r^*\CF$, see~\cite{br,dg}.

Recall that the geometric Satake equivalence takes an irreducible 
$G^\vee$-module $V^\blambda$ to the IC-sheaf $\IC(\ol\Gr{}^\blambda_G)$, and
the weight space $V^\blambda_\bmu$ is realized as $V^\blambda_\bmu=
\Phi_\bmu\IC(\ol\Gr{}^\blambda_G)=\Phi^\blambda_\bmu\IC(\ol\CW{}^\blambda_\bmu)$.
The usual stalks of both $\IC(\ol\Gr{}^\blambda_G)$ and 
$\IC(\ol\CW{}^\blambda_\bmu)$ at $\bmu$ are isomorphic up to shift to the 
associated graded $\on{gr}V^\blambda_\bmu$ with respect to the Brylinski-Kostant
filtration.

\subsection{The general case}
\label{general}
If we want to reconstruct the whole of $V^\blambda$ from the various slices
$\ol\CW{}_\bmu^\blambda$, we are missing the weight spaces $V^\blambda_\bmu$
with nondominant $\bmu$. To take care of the remaining weight spaces, for
arbitrary $\bmu$ we consider the moduli space $\ol\CW{}_\bmu^\blambda$ of 
the following data:

\textup{(a)} A $G$-bundle $\CP$ on $\BP^1$.

\textup{(b)} A trivialization $\sigma\colon \CP_{\on{triv}}|_{\BP^1\setminus\{0\}}
\iso\CP|_{\BP^1\setminus\{0\}}$ having a pole of degree $\leq\blambda$ at $0\in\BP^1$
(that is defining a point of $\ol\Gr{}_G^\blambda$).

\textup{(c)} A $B$-structure $\phi$ on $\CP$ of degree $w_0\bmu$ with the
fiber $B_-\subset G$ at $\infty\in\BP^1$ (with respect to the trivialization
$\sigma$ of $\CP$ at $\infty\in\BP^1$). Here $G\supset B_-\supset T$ is the
Borel subgroup opposite to $B$, and $w_0\in W_\fin$ is the longest element.

This construction goes back to~\cite{fm}. The space $\ol\CW{}_\bmu^\blambda$ is 
nonempty iff $\bmu\leq\blambda$. In this case it is
an irreducible affine normal Cohen-Macaulay variety of dimension
$\langle2\bar\rho^{\!\scriptscriptstyle\vee},\blambda-\bmu\rangle$, see~\cite{bfn3}.
In case $\bmu$ is dominant, the two definitions of $\ol\CW{}_\bmu^\blambda$
agree. At the other extreme, if $\blambda=0$, then $\ol\CW{}_{-\alpha}^0$ is
nothing but the open zastava space $\oZ^{-w_0\alpha}$.
The $T$-fixed point set $(\ol\CW{}_\bmu^\blambda)^T$ is nonempty iff the weight
space $V^\blambda_\bmu$ is not 0; in this case $(\ol\CW{}_\bmu^\blambda)^T$
consists of a single point denoted $\bmu$. We consider the repellent
$R^\blambda_\bmu\subset\ol\CW{}_\bmu^\blambda$. It is a closed subvariety of
dimension $\langle\bar\rho^{\!\scriptscriptstyle\vee},\blambda-\bmu\rangle$
(equidimensional). We have $V^\blambda_\bmu=\Phi_\bmu\IC(\ol\Gr{}^\blambda_G)=
\Phi^\blambda_\bmu\IC(\ol\CW{}^\blambda_\bmu)$, so that 
$V^\blambda=\bigoplus_{\bmu\in\Lambda}\Phi^\blambda_\bmu\IC(\ol\CW{}^\blambda_\bmu)$
(see~\cite{kr}).
Similarly to~\cite{bg}, one can introduce a crystal structure on the set
of irreducible components
$\bigsqcup_{\bmu\in\Lambda}\on{Irr}R^\blambda_\bmu$ (see~\cite{kr}),
so that the resulting crystal
is isomorphic to the integrable crystal ${\mathbf B}(\blambda)$
(for a beautiful survey on crystals, see~\cite{ka2}).

\subsection{Beilinson-Drinfeld slices}
\label{beilinson}
Let $\ul\blambda=(\blambda_1,\ldots,\blambda_N)$ be a collection of dominant
coweights of $G$. We consider the moduli space $\ul{\ol\CW}{}_\bmu^{\ul\blambda}$ of 
the following data:

\textup{(a)} A collection of points $(z_1,\ldots,z_N)\in\BA^N$ on the 
affine line $\BA^1\subset\BP^1$.

\textup{(b)} A $G$-bundle $\CP$ on $\BP^1$.

\textup{(c)} A trivialization 
$\sigma\colon \CP_{\on{triv}}|_{\BP^1\setminus\{z_1,\ldots,z_N\}}
\iso\CP|_{\BP^1\setminus\{z_1,\ldots,z_N\}}$ 
with a pole of degree $\leq\sum_{s=1}^N\blambda_s\cdot z_s$ on the complement.

\textup{(d)} A $B$-structure $\phi$ on $\CP$ of degree $w_0\bmu$ with the
fiber $B_-\subset G$ at $\infty\in\BP^1$ (with respect to the trivialization
$\sigma$ of $\CP$ at $\infty\in\BP^1$). 

$\ul{\ol\CW}{}_\bmu^{\ul\blambda}$ is nonempty iff 
$\bmu\leq\blambda:=\sum_{s=1}^N\blambda_s$. In this case it is an irreducible 
affine normal Cohen-Macaulay variety flat over $\BA^N$ of relative dimension 
$\langle2\bar\rho^{\!\scriptscriptstyle\vee},\blambda-\bmu\rangle$, see~\cite{bfn3}.
The fiber over $N\cdot0\in\BA^N$ is nothing but $\ol\CW{}_\bmu^\blambda$.
We can consider the Verdier specialization 
$\on{Sp}\IC(\ul{\ol\CW}{}_\bmu^{\ul\blambda})$ to the special fiber
$\ol\CW{}_\bmu^\blambda$. It is a perverse sheaf on 
$\ol\CW{}_\bmu^\blambda\times\BA^N$ smooth along the diagonal stratification of
$\BA^N$. We denote by $\Psi\IC(\ul{\ol\CW}{}_\bmu^{\ul\blambda})$ its restriction to
$\ol\CW{}_\bmu^\blambda\times\ul{z}$ where $\ul{z}$ is a point of $\BA^N_\BR$
such that $z_1>\ldots>z_N$. Then $$\Psi\IC(\ul{\ol\CW}{}_\bmu^{\ul\blambda})\simeq
\bigoplus_{\bmu\leq\bar{\nu}\leq\blambda,\ \bar{\nu}\in\Lambda^+}M^{\ul\blambda}_{\bar\nu}\otimes
\IC(\ol\CW{}^{\bar\nu}_\bmu),$$ where $M^{\ul\blambda}_{\bar\nu}$ is the multiplicity
$\on{Hom}_{G^\vee}(V^{\bar\nu},V^{\blambda_1}\otimes\ldots\otimes V^{\blambda_N})$.

\subsection{Convolution diagram over slices}
\label{convolution}
In the setup of~Section~\ref{beilinson} we consider the moduli space
$\ul{\wt\CW}{}_\bmu^{\ul\blambda}$ of the following data:

\textup{(a)} A collection of points $(z_1,\ldots,z_N)\in\BA^N$ on the 
affine line $\BA^1\subset\BP^1$.

\textup{(b)} A collection of $G$-bundles $(\CP_1,\ldots,\CP_N)$ on $\BP^1$.

\textup{(c)} A collection of isomorphisms  
$\sigma_s\colon \CP_{s-1}|_{\BP^1\setminus\{z_s\}}
\iso\CP_s|_{\BP^1\setminus\{z_s\}}$ 
with a pole of degree $\leq\blambda_s$ at $z_s$. Here $1\leq s\leq N$,
and $\CP_0:=\CP_{\on{triv}}$.

\textup{(d)} A $B$-structure $\phi$ on $\CP_N$ of degree $w_0\bmu$ with the
fiber $B_-\subset G$ at $\infty\in\BP^1$ (with respect to the trivialization
$\sigma_N\circ\ldots\circ\sigma_1$ of $\CP_N$ at $\infty\in\BP^1$). 

A natural projection 
$\varpi\colon \ul{\wt\CW}{}_\bmu^{\ul\blambda}\to\ul{\ol\CW}{}_\bmu^{\ul\blambda}$ 
sends $(\CP_1,\ldots,\CP_N,\sigma_1,\ldots,\sigma_N)$ to 
$(\CP_N,\sigma_N\circ\ldots\circ\sigma_1)$. We denote
$\varpi^{-1}(\ol\CW{}_\bmu^\blambda)$ by $\wt\CW{}_\bmu^{\ul\blambda}$. Then
$\varpi\colon \wt\CW{}_\bmu^{\ul\blambda}\to\ol\CW{}_\bmu^\blambda$ 
is stratified semismall, and $$\varpi_*\IC(\wt\CW{}_\bmu^{\ul\blambda})=
\bigoplus_{\bmu\leq\bar{\nu}\leq\blambda,\ \bar{\nu}\in\Lambda^+}
M^{\ul\blambda}_{\bar\nu}\otimes\IC(\ol\CW{}^{\bar\nu}_\bmu).$$

\section{Double affine Grassmannian}
\label{dag}
In this section $G$ is assumed to be a simply connected almost simple complex
algebraic group.

\subsection{The affine group and its Langlands dual}
We consider the minimal integral even positive definite $W_\fin$-invariant
symmetric bilinear form $(\cdot,\cdot)$ on the coweight lattice $\Lambda$. 
It gives rise to a central 
extension $\wh{G}$ of the polynomial version $G_{\BC[t^{\pm1}]}$ of the loop group:
$$1\to\BC^\times\to\wh{G}\to G_{\BC[t^{\pm1}]}\to1.$$
The loop rotation group $\BC^\times$ acts naturally on $G_{\BC[t^{\pm1}]}$, and this
action lifts to $\wh{G}$. We denote the corresponding semidirect product
$\BC^\times\ltimes\wh{G}$ by $G_\aff$. It is an untwisted affine Kac-Moody group
ind-scheme. 

We denote by $G_\aff^\vee$ the corresponding Langlands dual group.
Note that if $G$ is not simply laced, then $G_\aff^\vee$ is a twisted affine
Kac-Moody group, not to be confused with $(G^\vee)_\aff$. However, we have a
canonical embedding $G^\vee\hookrightarrow G_\aff^\vee$.

We fix a Cartan torus $\BC^\times\times T\times\BC^\times\subset G_\aff$ and its
dual Cartan torus $\BC^\times\times T^\vee\times\BC^\times\subset G_\aff^\vee$.
Here the first copy of $\BC^\times$ is the central $\BC^\times$, while the
second copy is the loop rotation $\BC^\times$. Accordingly, the weight lattice
$\Lambda_\aff$ of $G_\aff^\vee$ is $\BZ\oplus\Lambda\oplus\BZ$: the first copy
of $\BZ$ is the central charge (level), and the second copy is the energy.
A typical element $\lambda\in\Lambda_\aff$ will be written as
$\lambda=(k,\blambda,n)$.
The subset of dominant weights $\Lambda_\aff^+\subset\Lambda_\aff$ consists
of all the triples $(k,\blambda,n)$ such that $\blambda\in\Lambda^+$ and
$\langle\blambda,\bar\theta^{\scriptscriptstyle\vee}\rangle\leq k$. Here
$\bar\theta^{\scriptscriptstyle\vee}=\sum_{i\in I}a_i\alpha_i^{\scriptscriptstyle\vee}$ 
is the highest root of $G\supset B\supset T$.
We denote by $\Lambda^+_{\aff,k}\subset\Lambda^+_\aff$ the finite subset of 
dominant weights of level $k$; we also denote by 
$\Lambda_{\aff,k}\subset\Lambda_\aff$ the subset of all the weights of level $k$.
We say that $\lambda\geq\mu$ if $\lambda-\mu$ is an element of the submonoid
generated by the positive roots of $G_\aff^\vee$ (in particular, $\lambda$ and
$\mu$ must have the same level). Finally, let $\bar\omega_i,\ i\in I$, be
the fundamental coweights of $G$, and 
$\rho:=(1,0,0)+\sum_{i\in I}(a_i,\bar\omega_i,0)\in\Lambda_\aff$.

The affine Weyl group $W_\aff$ is the semidirect product $W_\fin\ltimes\Lambda$. 
For $k\in\BZ_{>0}$, we also consider its
version $W_{\aff,k}=W_\fin\ltimes k\Lambda$; it acts naturally on 
$\Lambda_{\aff,k}=\{k\}\times\Lambda\oplus\BZ$ (trivially on $\BZ$).
Every $W_{\aff,k}$-orbit on $\Lambda_{\aff,k}$ contains a unique representative
in $\Lambda_{\aff,k}^+$. It follows that if we denote by $\Gamma_k$ the group
of roots of unity of order $k$, then there is a natural isomorphism
$\Lambda_{\aff,k}^+/\BZ=W_{\aff,k}\backslash\Lambda\iso
\Hom(\Gamma_k,G)/\!\on{Ad}_G$.

\subsection{The quest}
\label{dream}
We would like to have a double affine Grassmannian $\Gr_{G_\aff}$ and a 
geometric Satake equivalence between the category of integrable representations
$\Rep(G_\aff^\vee)$ and an appropriate category of perverse sheaves on 
$\Gr_{G_\aff}$. Note that the affine Satake isomorphism at the level of functions
is established in~\cite{bk,bkp} (and in~\cite{garo} for arbitrary Kac-Moody
groups).

Such a quest was formulated by I.~Grojnowski in his talk at ICM-2006 in Madrid.
At approximately the same time, I.~Frenkel suggested that the integrable 
representations of level $k$ should be realized in cohomology of certain
instanton moduli spaces on $\BA^2/\Gamma_k$. Here $\Gamma_k$ acts on $\BA^2$
in a hyperbolic way: $\zeta(x,y)=(\zeta x,\zeta^{-1}y)$.

Note that the set of dominant coweights $\Lambda^+$ is well ordered, which 
reflects the fact that the affine Grassmannian $\Gr_G$ is an ind-projective
scheme. However, the set of affine dominant coweights $\Lambda_{\aff}^+$ is not 
well ordered: it does not have a minimal element. In fact, it has an 
automorphism group $\BZ$ acting by the energy shifts: 
$(k,\blambda,n)\mapsto(k,\blambda,n+n')$ (we add a multiple of the minimal
imaginary coroot $\delta$). This indicates that the sought for double affine
Grassmannian $\Gr_{G_\aff}$ is an object of semiinfinite nature.

At the moment, the only technical possibility of dealing with semiinfinite
spaces is via transversal slices to strata. Following I.~Frenkel's suggestion,
in the series~\cite{bf1,bf2,bf3} we developed a partial affine analogue of 
slices of~Section~\ref{slices} defined in terms of Uhlenbeck spaces 
$\CU_G(\BA^2/\Gamma_k)$.

\subsection{Dominant slices via Uhlenbeck spaces}
\label{domina}
The Uhlenbeck space $\CU_G^d(\BA^2)$ is a partial closure of the moduli space 
$\Bun^d_G(\BA^2)$ of $G$-bundles of second Chern class $d$ on the projective 
plane $\BP^2$ trivialized at the infinite line $\BP^1_\infty\subset\BP^2$, 
see~\cite{bfg}. It is known that $\Bun^d_G(\BA^2)$ is smooth quasiaffine,
and $\CU_G^d(\BA^2)$ is a connected affine variety of dimension 
$2dh_G^{\!\scriptscriptstyle\vee}$ (where $h_G^{\!\scriptscriptstyle\vee}$ is the dual
Coxeter number of $G$). Conjecturally, $\CU_G^d(\BA^2)$ is normal; in this 
case $\CU_G^d(\BA^2)$ is the affinization of $\Bun^d_G(\BA^2)$.

The group $G\times\GL(2)$ acts naturally on $\CU_G^d(\BA^2)$: the first
factor via the change of trivialization at $\BP^1_\infty$, and the second
factor via its action on $(\BP^2,\BP^1_\infty)$. The group $\Gamma_k$ is
embedded into $\GL(2)$. Given $\mu=(k,\bmu,m)\in\Lambda^+_{\aff,k}$ we choose its 
lift to a homomorphism from $\Gamma_k$ to $G$; thus $\Gamma_k$ embeds diagonally
into $G\times\GL(2)$ and acts on $\Bun^d_G(\BA^2)$. The fixed point subvariety
$\Bun^d_G(\BA^2)^{\Gamma_k}$ consists of $\Gamma_k$-equivariant bundles and
is denoted $\Bun^d_{G,\mu}(\BA^2/\Gamma_k)$; another choice of lift above leads 
to an isomorphic subvariety. Since $0\in\BA^2$ is a $\Gamma_k$-fixed point,
for any $\Gamma_k$-equivariant $G$-bundle $\CP\in\Bun^d_{G,\mu}(\BA^2/\Gamma_k)$
the group $\Gamma_k$ acts on the fiber $\CP_0$. This action defines an
element of $\Hom(\Gamma_k,G)/\!\on{Ad}_G$ to be denoted $[\CP_0]$.

Now given $\lambda=(k,\blambda,l)\in\Lambda^+_{\aff,k}$ we define
$\Bun^\lambda_{G,\mu}(\BA^2/\Gamma_k)$ as the subvariety of
$\Bun^d_{G,\mu}(\BA^2/\Gamma_k)$ formed by all $\CP$ such that the class
$[\CP_0]\in\Hom(\Gamma_k,G)/\!\on{Ad}_G$ is the image of $\lambda$, and
$d=k(l-m)+\frac{(\blambda,\blambda)-(\bmu,\bmu)}{2}$. It is a union of
connected components of $\Bun^d_{G,\mu}(\BA^2/\Gamma_k)$. Conjecturally,
$\Bun^\lambda_{G,\mu}(\BA^2/\Gamma_k)$ is connected. This conjecture is proved
if $G=\on{SL}(N)$, or $k=1$, or $k$ is big enough for arbitrary $G$ and 
fixed $\blambda,\bmu$.

Finally, we define the dominant slice $\ol\CW{}^\lambda_\mu$ as the closure 
$\CU^\lambda_{G,\mu}(\BA^2/\Gamma_k)$ of $\Bun^\lambda_{G,\mu}(\BA^2/\Gamma_k)$ in
the Uhlenbeck space $\CU_G^d(\BA^2)$.

\subsection{(Hyperbolic) stalks}
\label{hyperbolic}
The Cartan torus $T_\aff=\BC^\times\times T\times\BC^\times$ maps into 
$G\times\GL(2)$. Here the first copy of $\BC^\times$ goes to the diagonal torus 
of $\on{SL}(2)\subset\GL(2)$, while the second copy of $\BC^\times$ goes to the
center of $\GL(2)$. So $T_\aff$ acts on $\ol\CW{}^\lambda_\mu$, and we denote by
$\mu\in\ol\CW{}^\lambda_\mu$ the only fixed point. The corresponding repellent
$R^\lambda_\mu$ is the closed affine subvariety formed by all the points that flow
into $\mu$ under the action of $2\rho(t)$, as $t$ goes to $\infty$.
The corresponding hyperbolic stalk $\Phi^\lambda_\mu\IC(\ol\CW{}^\lambda_\mu)$ is 
conjecturally isomorphic to the weight space $V^\lambda_\mu$ of the integrable 
$G_\aff^\vee$-module $V^\lambda$ with highest weight $\lambda$. 
In type $A$ this conjecture
follows from the identification of $\ol\CW{}^\lambda_\mu$ with a Nakajima cyclic
quiver variety and I.~Frenkel's level-rank duality between the weight
multiplicities and the tensor product multiplicities~\cite{Fr,na1,na2,bf1}.
In type $ADE$ at level $1$ this conjecture follows from~\cite{bfn1}.
Also, as the notation suggests, the hyperbolic stalk
$\Phi^\lambda_\mu\IC(\ol\CW{}^\lambda_\mu)$ is isomorphic to the vanishing cycles
of $\IC(\ol\CW{}^\lambda_\mu)$ at $\mu$ with respect to a general function
vanishing at $\mu$~\cite{fk}.
The usual stalk of $\IC(\ol\CW{}^\lambda_\mu)$ at $\mu$ is conjecturally 
isomorphic to the associated graded of $V^\lambda_\mu$ with respect to the
the affine Brylinski-Kostant filtration~\cite{sl}. At level $1$, this conjecture
follows from the computation of the IC-stalks of Uhlenbeck spaces in~\cite{bfg}.

The affine analogs of generalized slices of 
Sections~\ref{general},~\ref{beilinson},~\ref{convolution} were constructed
in type $A$ in~\cite{bf2,bf3} in terms of Nakajima cyclic quiver
varieties mentioned above. For arbitrary $G$, the desired generalized slices
are expected to be the Uhlenbeck partial compactifications of the moduli spaces 
of $\Gamma_k$-equivariant $G_c$-instantons (where $G_c\subset G$ is a maximal
compact subgroup) on multi Taub-NUT spaces (for a physical explanation
via a supersymmetric conformal field theory in $6$ dimensions, see~\cite{wi}).
Unfortunately, we are still lacking a modular definition of the Uhlenbeck
compactification~\cite{ba}, and the existing {\em ad hoc} constructions are
not flexible enough. Another approach via the Coulomb branches of framed affine
quiver gauge theories following~\cite{na3,bfn2,bfn3,bfn4} is described in the
remaining sections. For a beautiful short introduction to the Coulomb branches,
the reader may consult~\cite{na5,na4}.

\section{Coulomb branches of $3d\ \CN=4$ quiver gauge theories}
\label{coulomb}

\subsection{General setup}
\label{general setup}
Let $\bN$ be a finite dimensional representation of a complex connected
reductive group $\sG$ (having nothing to do with $G$ of previous sections).  
We consider the moduli space $\CR_{\sG,\bN}$ of triples $(\CP,\sigma,s)$
where $\CP$ is a $\sG$-bundle on the formal disc $D=\on{Spec}\CO;\ \sigma$
is a trivialization of $\CP$ on the punctured formal disc $D^*=\on{Spec}\CK$;
and $s$ is a section of the associated vector bundle 
$\CP_{\on{triv}}\Gtimes\bN$ on $D^*$ such that $s$ extends to a 
regular section of $\CP_{\on{triv}}\Gtimes\bN$ on $D$, and
$\sigma(s)$ extends to a regular section of $\CP\Gtimes\bN$ on 
$D$. In other words, $s$ extends to a regular section of the vector bundle
associated to the $\sG$-bundle
glued from $\CP$ and $\CP_{\on{triv}}$ on the non-separated formal scheme glued
from $2$ copies of $D$ along $D^*$ ({\em raviolo}).
The group $\sG_\CO$ acts on $\CR_{\sG,\bN}$ by changing the trivialization 
$\sigma$, and we have an evident 
$\sG_\CO$-equivariant projection $\CR_{\sG,\bN}\to\Gr_\sG$ forgetting $s$.
The fibers of this projection are profinite dimensional vector spaces:
the fiber over the base point is $\bN\otimes\CO$, and all the other fibers
are subspaces in $\bN\otimes\CO$ of finite codimension. One may say that
$\CR_{\sG,\bN}$ is a $\sG_\CO$-equvariant ``constructible profinite dimensional
vector bundle'' over $\Gr_\sG$. The $\sG_\CO$-equivariant Borel-Moore homology
$H^{\sG_\CO}_\bullet(\CR_{\sG,\bN})$ is well-defined, and forms an associative 
algebra with respect to a convolution operation. This algebra is commutative, 
finitely generated and integral, and its spectrum 
$\CM_C(\sG,\bN)=\on{Spec}H^{\sG_\CO}_\bullet(\CR_{\sG,\bN})$ is an irreducible
normal affine variety of dimension $2\on{rk}(\sG)$, the {\em Coulomb branch}.
It is supposed to be a (singular) hyper-K\"ahler manifold~\cite{sw}.

Let $\sT\subset\sG$ be a Cartan torus with Lie algebra $\ft\subset\fg$.
Let $\sW=N_\sG(\sT)/\sT$ be the corresponding Weyl group. Then the equivariant
cohomology $H^\bullet_{\sG_\CO}(\on{pt})=\BC[\ft/\sW]$ forms a subalgebra of
$H^{\sG_\CO}_\bullet(\CR_{\sG,\bN})$ (a {\em Cartan subalgebra}), so we have a 
projection $\varPi\colon \CM_C(\sG,\bN)\to\ft/\sW$.

Finally, the algebra $H^{\sG_\CO}_\bullet(\CR_{\sG,\bN})$ comes equipped with 
quantization: a $\BC[\hbar]$-deformation 
$\BC_\hbar[\CM_C(\sG,\bN)]=H^{\BC^\times\ltimes\sG_\CO}_\bullet(\CR_{\sG,\bN})$
where $\BC^\times$ acts by loop rotations, and 
$\BC[\hbar]=H^\bullet_{\BC^\times}(\on{pt})$. It gives rise to a Poisson bracket
on $\BC[\CM_C(\sG,\bN)]$ with an open symplectic leaf, so that $\varPi$ 
becomes an integrable system: $\BC[\ft/\sW]\subset\BC[\CM_C(\sG,\bN)]$ is a 
Poisson-commutative polynomial subalgebra with $\on{rk}(\sG)$ generators.

\subsection{Flavor symmetry}
\label{flavor}
Suppose we have an extension $1\to\sG\to\tilde\sG\to\sG_F\to1$ where $\sG_F$
is a connected reductive group (a {\em flavor group}), and the action of
$\sG$ on $\bN$ is extended to an action of $\tilde\sG$. Then the action of
$\sG_\CO$ on $\CR_{G,\bN}$ extends to an action of $\tilde\sG_\CO$, and the
convolution product defines a commutative algebra structure on the equivariant
Borel-Moore homology $H^{\tilde\sG_\CO}_\bullet(\CR_{\sG,\bN})$. We have the
restriction homomorphism $H^{\tilde\sG_\CO}_\bullet(\CR_{\sG,\bN})\to 
H^{\sG_\CO}_\bullet(\CR_{\sG,\bN})=H^{\tilde\sG_\CO}_\bullet(\CR_{\sG,\bN})
\otimes_{H^\bullet_{\sG_F}(\on{pt})}\BC$. In other words, $\ul\CM{}_C(\sG,\bN):=
\Spec H^{\tilde\sG_\CO}_\bullet(\CR_{\sG,\bN})$ is a deformation 
of $\CM_C(\sG,\bN)$ over $\Spec H^\bullet_{\sG_F}(\on{pt})=\ft_F/\sW_F$.

We will need the following version of this construction. Let $\sZ\subset\sG_F$
be a torus embedded into the flavor group. We denote by $\tilde\sG{}^\sZ$
the pullback extension $1\to\sG\to\tilde\sG{}^\sZ\to\sZ\to1$. We define
$\ul\CM{}_C^\sZ(\sG,\bN):=\Spec H^{\tilde\sG{}^\sZ_\CO}_\bullet(\CR_{\sG,\bN})$:
a deformation of $\CM_C(\sG,\bN)$ over $\fz:=\Spec H^\bullet_\sZ(\on{pt})$.

Since $\CM_C(\sG,\bN)$ is supposed to be a hyper-K\"ahler manifold, its
flavor deformation should come together with a (partial) resolution.
To construct it, we consider the obvious projection
$\tilde\pi\colon \CR_{\tilde\sG,\bN}\to\Gr_{\tilde\sG}\to\Gr_{\sG_F}$. 
Given a dominant coweight $\lambda_F\in\Lambda^+_F\subset\Gr_{\sG_F}$, we set
$\CR_{\tilde\sG,\bN}^{\lambda_F}:=\tilde\pi^{-1}(\lambda_F)$, and consider
the equivariant Borel-Moore homology
$H^{\tilde\sG{}^\sZ_\CO}_\bullet(\CR_{\tilde\sG,\bN}^{\lambda_F})$. It carries a convolution
module structure over $H^{\tilde\sG{}^\sZ_\CO}_\bullet(\CR_{\sG,\bN})$. We consider 
$\wt{\ul\CM}{}_C^{\sZ,\lambda_F}(\sG,\bN):=
\Proj(\bigoplus_{n\in\BN}H^{\tilde\sG{}^\sZ_\CO}_\bullet(\CR_{\tilde\sG,\bN}^{n\lambda_F}))
\stackrel{\varpi}{\longrightarrow}\ul\CM{}_C^\sZ(\sG,\bN)$.
We denote $\varpi^{-1}(\CM_C(\sG,\bN))$ by $\wt\CM_C^{\lambda_F}(\sG,\bN)$.
We have $\wt\CM_C^{\lambda_F}(\sG,\bN)=
\Proj(\bigoplus_{n\in\BN}H^{\sG_\CO}_\bullet(\CR_{\tilde\sG,\bN}^{n\lambda_F}))$.

More generally, for a strictly convex (i.e.\ not containing nontrivial 
subgroups) cone $\sV\subset\Lambda^+_F$, we consider the multi projective 
spectra $\wt{\ul\CM}{}_C^{\sZ,\sV}(\sG,\bN):=
\Proj(\bigoplus_{\lambda_F\in\sV}H^{\tilde\sG{}^\sZ_\CO}_\bullet(\CR_{\tilde\sG,\bN}^{\lambda_F}))
\stackrel{\varpi}{\longrightarrow}\ul\CM{}_C^\sZ(\sG,\bN)$ and
$\wt\CM_C^\sV(\sG,\bN):=
\Proj(\bigoplus_{\lambda_F\in\sV}H^{\sG_\CO}_\bullet(\CR_{\tilde\sG,\bN}^{\lambda_F}))
\stackrel{\varpi}{\longrightarrow}\CM_C(\sG,\bN)$.

\subsection{Quiver gauge theories}
\label{quiver}
Let $Q$ be a quiver with $Q_0$ the set of vertices, and $Q_1$ the set of arrows.
An arrow $e\in Q_1$ goes from its tail $t(e)\in Q_0$ to its head $h(e)\in Q_0$. 
We choose a 
$Q_0$-graded vector spaces $V:=\bigoplus_{j\in Q_0}V_j$ and 
$W:=\bigoplus_{j\in Q_0}W_j$. We set $\sG=\GL(V):=\prod_{j\in Q_0}\GL(V_j)$. 
We choose a second grading
$W=\bigoplus_{s=1}^NW^{(s)}$ compatible with the $Q_0$-grading of $W$.
We set $\sG_F$ to be a Levi subgroup $\prod_{s=1}^N\prod_{j\in Q_0}\GL(W^{(s)}_j)$
of $\GL(W)$, and $\tilde\sG:=\sG\times\sG_F$. Finally, we define a central
subgroup $\sZ\subset\sG_F$ as follows: $\sZ:=\prod_{s=1}^N\Delta_{\BC^\times}^{(s)}
\subset\prod_{s=1}^N\prod_{j\in Q_0}\GL(W^{(s)}_j)$, where
$\BC^\times\cong\Delta_{\BC^\times}^{(s)}\subset\prod_{j\in Q_0}\GL(W^{(s)}_j)$ is the
diagonally embedded subgroup of scalar matrices.
The reductive group $\tilde\sG$ acts naturally on $\bN:=
\bigoplus_{e\in Q_1}\Hom(V_{t(e)},V_{h(e)})\oplus\bigoplus_{j\in Q_0}\Hom(W_j,V_j)$. 

The Higgs branch of the corresponding quiver gauge theory is the Nakajima
quiver variety $\CM_H(\sG,\bN)=\fM(V,W)$. We are interested in the Coulomb 
branch $\CM_C(\sG,\bN)$.

\subsection{Back to slices in an affine Grassmannian}
\label{back}
Let now $G$ be an adjoint simple simply laced algebraic group. We choose
an orientation $\Omega$ of its Dynkin graph (of type $ADE$), and denote by
$I$ its set of vertices. Given an $I$-graded vector space $W$ we encode its
dimension by a dominant coweight  
$\blambda:=\sum_{i\in I}\dim(W_i)\bar{\omega}_i\in\Lambda^+$ of $G$. 
Given an $I$-graded vector space $V$ we encode its dimension by a positive
coroot combination $\alpha:=\sum_{i\in I}\dim(V_i)\alpha_i\in\Lambda_+$.
We set $\bmu:=\blambda-\alpha\in\Lambda$. Given a direct sum decomposition
$W=\bigoplus_{s=1}^NW^{(s)}$ compatible with the $I$-grading of $W$ as
in~Section~\ref{quiver}, we set 
$\blambda_s:=\sum_{i\in I}\dim(W_i^{(s)})\bar{\omega}_i\in\Lambda^+$, and
finally, $\ul\blambda:=(\blambda_1,\ldots,\blambda_N)$.

Recall the notations of~Section~\ref{flavor}. Since the flavor group 
$\sG_F$ is a Levi subgroup of $\GL(W)$, its weight lattice is naturally
identified with $\BZ^{\dim W}$. More precisely, we choose a basis
$w_1,\ldots,w_{\dim W}$ of $W$ such that any $W_i,\ i\in I$, and 
$W^{(s)},\ 1\leq s\leq N$, is spanned by a subset of the basis, and we assume
the following monotonicity condition: if for $1\leq a<b<c\leq\dim W$ we have
$w_a,w_b\in W^{(s)}$ for certain $s$, then $w_b\in W^{(s)}$ as well.
We define a strictly convex cone 
$\sV=\{(n_1,\ldots,n_{\dim W})\}\subset\Lambda^+_F\subset\BZ^{\dim W}$ 
by the following conditions: (a) if $w_k\in W^{(s)},\ w_l\in W^{(t)}$, 
and $s<t$, then $n_k\geq n_l\geq0$; (b) if $w_k,w_l\in W^{(s)}$, then $n_k=n_l$.
The following isomorphisms are constructed in~\cite{bfn3}
(notations of~Section~\ref{slices}):
$$\ol\CW{}_\bmu^\blambda\iso\CM_C(\sG,\bN),\
\ul{\ol\CW}{}_\bmu^{\ul\blambda}\iso\ul\CM{}_C^\sZ(\sG,\bN),$$
(we learned of their existence from V.~Pestun).
We also expect the following isomorphisms:
$$\ul{\wt\CW}{}_\bmu^{\ul\blambda}\iso\wt{\ul\CM}{}_C^{\sZ,\sV}(\sG,\bN),\
\wt\CW{}_\bmu^{\ul\blambda}\iso\wt\CM{}_C^\sV(\sG,\bN).$$

In case $G$ is an adjoint simple non simply laced algebraic group, it can be
obtained by folding from a simple simply laced group $\tilde G$ (i.e.\ as
the fixed point set of an outer automorphism of $\tilde G$). The corresponding
automorphism of the Dynkin quiver of $\tilde G$  acts on the above
Coulomb branches, and the slices for $G$ can be realized as the fixed point
sets of these Coulomb branches.

\subsection{Back to slices in a double affine Grassmannian}
\label{double back}
We choose an orientation of an {\em affine} Dynkin graph of type 
$A^{(1)},D^{(1)},E^{(1)}$ with the set of vertices $\tilde{I}=I\sqcup\{i_0\}$,
and repeat the construction of~Section~\ref{back} for an affine dominant
coweight 
$\lambda=\sum_{i\in\tilde{I}}\dim(W_i)\omega_i=(k,\blambda,0)\in\Lambda_\aff^+$,
a positive coroot combination 
$\alpha=\sum_{i\in\tilde{I}}\dim(V_i)\alpha_i\in\Lambda_{\aff,+}$, and
$\mu:=\lambda-\alpha=(k,\bmu,n)\in\Lambda_\aff$.

We define the slices in $\Gr_{G_\aff}$ (where $G$ is the corresponding adjoint
simple simply laced algebraic group) as
$$\ol\CW{}_\mu^\lambda:=\CM_C(\sG,\bN),\
\ul{\ol\CW}{}_\mu^{\ul\lambda}:=\ul\CM{}_C^\sZ(\sG,\bN),\
\ul{\wt\CW}{}_\mu^{\ul\lambda}:=\wt{\ul\CM}{}_C^{\sZ,\sV}(\sG,\bN),\
\wt\CW{}_\mu^{\ul\lambda}:=\wt\CM{}_C^\sV(\sG,\bN).$$
If $\mu$ is dominant, the slices $\ol\CW{}_\mu^\lambda$ conjecturally coincide
with the ones of~Section~\ref{domina}. In type $A$ this conjecture follows from
the computation~\cite{nt} of Coulomb branches of the cyclic quiver gauge 
theories and their identification with the Nakajima cyclic quiver varieties.

Note that $\pi_0(\CR_{\sG,\bN})=\pi_0(\Gr_{\GL(V)})=\pi_1(\GL(V))=\BZ^{\tilde{I}}$,
so that 
$H^{\sG_\CO}_\bullet(\CR_{\sG,\bN})=\BC[\CM_C(\sG,\bN)]=\BC[\ol\CW{}_\mu^\lambda]$
is $\BZ^{\tilde{I}}$-graded. We identify $\BZ^{\tilde{I}}$ with the root lattice of 
$T_\aff\subset G_\aff\colon \BZ^{\tilde{I}}=
\BZ\langle\alpha^{\!\scriptscriptstyle\vee}_i\rangle_{i\in\tilde{I}}$. Then the
$\BZ^{\tilde{I}}$-grading on $\BC[\ol\CW{}_\mu^\lambda]$ corresponds to a
$T_\aff$-action on $\ol\CW{}_\mu^\lambda$. Composing with the cocharacter
$2\rho\colon \BC^\times\to T_\aff$, we obtain an action of $\BC^\times$ on
$\ol\CW{}_\mu^\lambda$. Conjecturally, the fixed point set 
$(\ol\CW{}_\mu^\lambda)^{\BC^\times}$ is nonempty iff the $V^\lambda_\mu\ne0$, and 
in this case the fixed point set consists of a single point denoted by $\mu$.
We consider the corresponding repellent $R^\lambda_\mu\subset\ol\CW{}_\mu^\lambda$
and the hyperbolic stalk $\Phi^\lambda_\mu\IC(\ol\CW{}_\mu^\lambda)$.

Similarly to~Section~\ref{back}, in case $G$ is an adjoint simple non simply
laced group, the Dynkin diagram of its affinization can be obtained by folding
of a Dynkin graph of type $A^{(1)},D^{(1)},E^{(1)}$, and the above slices for
$G$ are defined as the fixed point sets of the corresponding slices for the
unfolding of $G$. The repellents and the hyperbolic stalks are thus defined 
for arbitrary simple $G$ too, and we expect the conclusions 
of~Sections~\ref{general},~\ref{beilinson},~\ref{convolution} to hold in the 
affine case as well.

\subsection{Warning}
\label{warning}
In order to formulate the statements about multiplicities for fusion and 
convolution as in~\ref{beilinson} and~\ref{convolution}, we must have
closed embeddings of slices
$\ol\CW{}_\mu^{\lambda'}\hookrightarrow\ol\CW{}_\mu^\lambda$ for 
$\lambda'\leq\lambda\in\Lambda^+_\aff$. Certainly we do have the natural
closed embeddings of generalized slices in $\Gr_G\colon
\ol\CW{}_\bmu^{\blambda'}\hookrightarrow\ol\CW{}_\bmu^\blambda,\ 
\blambda'\leq\blambda\in\Lambda^+$, but these embeddings have no manifest
interpretation in terms of Coulomb branches (see~Section~\ref{multiplication}
below for a partial advance, though). For a slice in $\Gr_G$, 
the collection of closures of symplectic leaves in $\ol\CW{}_\bmu^\blambda$
coincides with the collection of smaller slices
$\ol\CW{}_\bmu^{\blambda'}\subset\ol\CW{}_\bmu^\blambda,\ 
\bmu\leq\blambda'\leq\blambda,\ \blambda'\in\Lambda^+$.
However, in the affine case, in general there are {\em more} symplectic
leaves in $\ol\CW{}_\mu^\lambda$ than the cardinality of
$\{\lambda'\in\Lambda^+_\aff: \mu\leq\lambda'\leq\lambda\}$. For example,
if $k=1$, and $\mu=0$, so that $\ol\CW{}^\lambda_\mu\simeq\CU^d_G(\BA^2)$, the
symplectic leaves are numbered by the partitions of size $\leq d$: they are
all of the form $\CS\times\Bun_G^{d'}(\BA^2)$ where $0\leq d'\leq d$, and
$\CS$ is a stratum of the diagonal stratification of $\on{Sym}^{d-d'}\BA^2$.

Thus we expect that the slice $\ol\CW{}_\mu^{\lambda'}$ for
$\lambda'\in\Lambda^+_\aff,\ \mu\leq\lambda'\leq\lambda$, is isomorphic to
the closure of a symplectic leaf in $\ol\CW{}_\mu^\lambda$.
We also do expect the multiplicity of $\IC(\ol\CW{}^{\lambda'}_\mu)$ in
$\Psi\IC(\ul{\ol\CW}{}_\mu^{\ul\lambda})=\varpi_*\IC(\wt\CW{}_\mu^{\ul\lambda})$ 
to be $M^{\ul\lambda}_{\lambda'}=
\on{Hom}_{G_\aff^\vee}(V^{\lambda'},V^{\lambda_1}\otimes\ldots\otimes V^{\lambda_N})$
for any $\lambda'\in\Lambda^+_\aff$ such that $\mu\leq\lambda'\leq\lambda$.
However, it is possible 
that the IC sheaves of {\em other} symplectic leaves' closures also enter 
$\Psi\IC(\ul{\ol\CW}{}_\mu^{\ul\lambda})=\varpi_*\IC(\wt\CW{}_\mu^{\ul\lambda})$
with nonzero multiplicities. We should understand the representation-theoretic
meaning of these extra multiplicities, 
cf.~\cite[Theorem~5.15 and~Remark~5.17(3)]{na2} for $G$ of type $A$.

Also, the closed embeddings of slices (for Levi subgroups of
$G_\aff$) seem an indispensable tool for constructing a $\fg_\aff^\vee$-action
on $\bigoplus_\mu\Phi^\lambda_\mu\IC(\ol\CW{}_\mu^\lambda)$ or a structure of
$\fg_\aff^\vee$-crystal on $\bigsqcup_\mu\on{Irr}(R^\lambda_\mu)$ (via reduction
to Levi subgroups), cf.~\cite{kr}.

\subsection{Further problems}
\label{fusion}
Note that the construction of~Section~\ref{double back} uses no specific properties
of the affine Dynkin graphs, and works in the generality of arbitrary graph $Q$
without edge loops and the corresponding Kac-Moody Lie algebra $\fg_Q$. 
We still expect the
conclusions of~Sections~\ref{general},~\ref{beilinson},~\ref{convolution} to 
hold in this generality, see~\cite[3(x)]{bfn3}. 

The only specific feature of the affine case is as follows.
Recall that the category $\on{Rep}_k(G^\vee_\aff)$ of integrable 
$G^\vee_\aff$-modules at level $k\in\BZ_{>0}$ is equipped with a braided balanced
tensor {\em fusion} structure~\cite{ms,baki}. Unfortunately, I have no clue how
this structure is reflected in the geometry of $\Gr_{G_\aff}$. I believe this
is one of the most pressing problems about $\Gr_{G_\aff}$.


\section{Applications}
\label{applications}

\subsection{Hikita conjecture}
\label{hikita}
We already mentioned in~Section~\ref{double back} that in case $V^\lambda_\mu\ne0$
we expect the fixed point set $(\ol\CW{}^\lambda_\mu)^{T_\aff}$ to consist of a
single point $\mu$. This point is the support of a nilpotent scheme $\bfmu$
defined as follows: we choose a $T_\aff$-equivariant embedding
$\ol\CW{}^\lambda_\mu\hookrightarrow\BA^N$ into a representation of $T_\aff$,
and define $\bfmu$ as the scheme-theoretic intersection of $\ol\CW{}^\lambda_\mu$
with the zero weight subspace $\BA^N_0$ inside $\BA^N$. The resulting 
subscheme $\bfmu\subset\ol\CW{}^\lambda_\mu$ is independent of the choice of
a $T_\aff$-equivariant embedding $\ol\CW{}^\lambda_\mu\hookrightarrow\BA^N$.
According to the Hikita conjecture~\cite{hi}, the ring $\BC[\bfmu]$ is 
expected to be isomorphic to the cohomology ring $H^\bullet(\fM(V,W))$ of the
corresponding Nakajima affine quiver variety, see~Section~\ref{quiver}.
This is an instance of {\em symplectic duality} ($3d$ mirror symmetry) between 
Coulomb and Higgs branches. The Hikita conjecture for
the slices $\ol\CW{}^\blambda_\bmu$ in $\Gr_G$ and the corresponding finite type
Nakajima quiver varieties is proved in~\cite{ktwy} for types $A,D$
(and conditionally for types $E$).

\subsection{Monopole formula}
\label{monopole}
We return to the setup of~Section~\ref{general setup}. Recall that $\CR_{\sG,\bN}$
is a union of (profinite dimensional) vector bundles over $\sG_\CO$-orbits
in $\Gr_\sG$. The corresponding Cousin spectral sequence converging to
$H_\bullet^{\sG_\CO}(\CR_{\sG,\bN})$ degenerates and allows to compute the 
equivariant Poincar\'e polynomial (or rather Hilbert series) 
\begin{equation}
\label{mono}
P_t^{\sG_\CO}(\CR_{\sG,\bN})=
\sum_{\theta\in\Lambda^+_\sG}t^{d_\theta-2\langle\rho_\sG,\theta\rangle}P_\sG(t;\theta).
\end{equation}
Here $\deg(t)=2,\ P_\sG(t;\theta)=\prod(1-t^{d_i})^{-1}$ is the Hilbert series
of the equivariant cohomology $H^\bullet_{\on{Stab}_\sG(\theta)}(\on{pt})\ (d_i$
are the degrees of generators of the ring of $\on{Stab}_\sG(\theta)$-invariant 
functions on its Lie algebra), and $d_\theta=
\sum_{\chi\in\Lambda^\vee_\sG}\on{max}(-\langle\chi,\theta\rangle,0)\dim\bN_\chi$.
This is a slight variation of the {\em monopole formula} of~\cite{chz}.
Note that the series~(\ref{mono}) may well diverge (even as a formal Laurent
series: the space of homology of given degree may be infinite-dimensional),
e.g.\ this is always the case for unframed quiver gauge theories.
To ensure its convergence (as a formal Taylor series with the constant term~1)
one has to impose the so called `good' or `ugly' assumption on the theory.
In this case the resulting $\BN$-grading on $H_\bullet^{\sG_\CO}(\CR_{\sG,\bN})$
gives rise to a $\BC^\times$-action on $\CM_C(\sG,\bN)$, making it a conical
variety with a single (attracting) fixed point.

Now recall the setup of~Sections~\ref{quiver},~\ref{back}; in particular,
the isomorphism $\ol\CW{}_\bmu^\blambda\iso\CM_C(\sG,\bN)$. In case $\bmu$ is
dominant, the slice $\ol\CW{}_\bmu^\blambda\subset\Gr_G$ is conical with
respect to the loop rotation $\BC^\times$-action. However, this action is
{\em not} the one of the previous paragraph. They differ by a hamiltonian
$\BC^\times$-action (preserving the Poisson structure). The Hilbert series of
$\ol\CW{}_\bmu^\blambda$ graded by the loop rotation $\BC^\times$-action is given by
\begin{equation}
\label{loop}
P_t(\BC[\ol\CW{}_\bmu^\blambda])=\sum_{\theta\in\Lambda^+_\sG}
t^{d_\theta-2\langle\rho_\sG,\theta\rangle-\frac{1}{2}\bar\theta^\dagger\cdot\det\bN_{\on{hor}}+
\frac{1}{2}\bar\theta^\dagger\cdot C\cdot\alpha}P_\sG(t;\theta).
\end{equation}
Here $\deg(t)=1;\ \alpha=\blambda-\bmu\in\Lambda_+=\BN^I;\ \bar\theta$ is the 
class of $\theta\in\Lambda_\sG=\Lambda_{\GL(V)}$ in $\pi_0\Gr_{\GL(V)}=\BZ^I;\ 
\bar\theta^\dagger$ is the transposed row-vector; $C$ is the $I\times I$ Cartan
matrix of $G$; and $\bN_{\on{hor}}=\bigoplus_{i\to j\in\Omega}\Hom(V_i,V_j)$ is the 
``horizontal'' summand of $\GL(V)$-module $\bN$, so that $\det\bN_{\on{hor}}$
is a character of $\GL(V)$, i.e.\ an element of $\BZ^I$.

Finally, we consider a double affine Grassmannian slice $\ol\CW{}_\mu^\lambda$
with dominant $\mu$ as in~Section~\ref{domina}. The analogue of the loop
rotation action of the previous paragraph is the action of the second copy
of $\BC^\times$ (the center of $\GL(2)$) in~Section~\ref{hyperbolic}. We expect
that the Hilbert series of $\ol\CW{}_\mu^\lambda$ graded by this 
$\BC^\times$-action is given by the evident affine analogue of the 
formula~(\ref{loop}) (with the $\tilde{I}\times\tilde{I}$ Cartan matrix
$C_\aff$ of $G_\aff$ replacing $C$). In particular, in case of level 1, this 
gives a formula for the Hilbert series of the coordinate ring
$\BC[\CU^d_G(\BA^2)]$ of the Uhlenbeck space proposed in~\cite{cfhm}.
Note that the latter formula works for arbitrary $G$, not necessarily simply 
laced one. In type $A$ it follows from the results of~\cite{nt}.

\subsection{Zastava}
\label{zastava}
Let us consider the Coulomb branch $\CM_C(\sG,\bN)$ of an {\em unframed}
quiver gauge theory for an $ADE$ type quiver: $W_i=0\ \forall i\in I$, 
so that $\bN=\bN_{\on{hor}}$.
An isomorphism $\CM_C(\sG,\bN)\iso\oZ^\alpha$ with the open 
zastava\footnote{Zastava $=$ {\em flags} in Croatian.}
(the moduli space of degree $\alpha$ based maps from the projective line
$\BP^1\ni\infty$ to the flag variety $\CB\ni B_-$ of $G$, where
$\alpha=\sum_{i\in I}(\dim V_i)\alpha_i$), is constructed in~\cite{bfn3}
(we learned of its existence from V.~Pestun). As the name suggests, the open 
zastava is a (dense smooth symplectic) open subvariety in the zastava space 
$Z^\alpha$, a normal Cohen-Macaulay affine Poisson variety. 

Note that there is another version of zastava ${\mathbf Z}^\alpha$ that is
the solution of a moduli problem ($G$-bundles on $\BP^1$ with a generalized
$B$-structure and an extra $U_-$-structure transversal 
at $\infty\in\BP^1$)~\cite{bfgm} given by a scheme cut out by the Pl\"ucker
equations. This scheme is not reduced in general (the first example occurs
in type $A_4$)~\cite{fema}, and $Z^\alpha$ is the corresponding variety:
$Z^\alpha:={\mathbf Z}^\alpha_{\on{red}}$.

We already mentioned that the open zastava is a particular case of a 
generalized slice: $\oZ^\alpha=\ol\CW{}^0_{w_0\alpha}$. The zastava space $Z^\alpha$ 
is the limit of slices in the following sense: for any $\blambda\geq\bmu$ such 
that $w_0\bmu-w_0\blambda=\alpha$, there is a loop rotation equivariant regular
birational morphism 
$s^\blambda_\bmu\colon \ol\CW{}^{-w_0\blambda}_{-w_0\bmu}\to Z^\alpha$, and for any
$N\in\BN$ and big enough dominant $\bmu$, the corresponding morphism of the
coordinate rings graded by the loop rotations 
$(s^\blambda_\bmu)^*\colon \BC[Z^\alpha]\to\BC[\ol\CW{}^{-w_0\blambda}_{-w_0\bmu}]$
is an isomorphism in degrees $\leq N$ (both $\BC[Z^\alpha]$ and 
$\BC[\ol\CW{}^{-w_0\blambda}_{-w_0\bmu}]$ for dominant $\bmu$ are positively graded).

Now $\BC[Z^\alpha]$ is obtained by the following version of the Coulomb branch
construction. Given a vector space $U$ we define the positive part of the 
affine Grassmannian $\Gr^+_{GL(U)}\subset\Gr_{GL(U)}$ as the moduli space of 
vector bundles $\CU$ on the formal disc $D=\on{Spec}(\CO)$ equipped with 
trivialization $\sigma\colon \CU|_{D^*}\iso U\otimes\CO_{D^*}$ on the formal
punctured disc $D^*=\on{Spec}(\CK)$ such that $\sigma$ extends through the
puncture as an embedding $\sigma\colon \CU\hookrightarrow U\otimes\CO_D$.
Since $\sG=\on{GL}(V)=\prod_{i\in I}\on{GL}(V_i)$, we have 
$\Gr_{\on{GL}(V)}=\prod_{i\in I}\Gr_{\on{GL}(V_i)}$, and we define 
$\Gr^+_{\on{GL}(V)}=\prod_{i\in I}\Gr^+_{\on{GL}(V_i)}$. Finally, we define
$\CR^+_{\sG,\bN}$ as the preimage of $\Gr^+_{\on{GL}(V)}\subset\Gr_{\on{GL}(V)}$ under
$\CR_{\sG,\bN}\to\Gr_{\on{GL}(V)}$. Then $H_\bullet^{\sG_\CO}(\CR^+_{\sG,\bN})$ forms a
convolution subalgebra of $H_\bullet^{\sG_\CO}(\CR_{\sG,\bN})$, and an isomorphism
$\CM_C^+(\sG,\bN):=\on{Spec}H_\bullet^{\sG_\CO}(\CR^+_{\sG,\bN})\iso Z^\alpha$
is constructed in~\cite{bfn3}. 

An analogue of the monopole formula~(\ref{loop}) gives the character of the
$T\times\BC^\times$-module $\BC[Z^\alpha]$:
\begin{equation}
\label{cartan loop}
\chi(\BC[Z^\alpha])=\sum_{\Lambda_\sG^{++}}z^{\bar\theta}
t^{d_\theta-2\langle\rho_\sG,\theta\rangle-\frac{1}{2}\bar\theta^\dagger\cdot\det\bN+
\frac{1}{2}\bar\theta^\dagger\cdot C\cdot\alpha}P_\sG(t;\theta).
\end{equation}
Here $\Lambda_\sG^{++}$ is the set of $I$-tuples of partitions; $i$-th partition
having length at most $\dim V_i$ (recall that the cone of dominant coweights
$\Lambda_\sG^+$ is formed by the $I$-tuples of nonincreasing sequences 
$(\lambda^{(i)}_1\geq\lambda^{(i)}_2\geq\ldots\geq\lambda^{(i)}_{\dim V_i})$ of
integers, and for $\Lambda_\sG^{++}\subset\Lambda_\sG^+$ we require these integers
to be nonnegative). Also, $z$ denotes the coordinates on the Cartan torus
$T\subset G$ identified with $(\BC^\times)^I$ via 
$z_i=\alpha^{\scriptscriptstyle\vee}_i$.

The character of the $T\times\BC^\times$-module $\BC[Z^\alpha]$ for $G$
of type $ADE$ was also computed in~\cite{bf4}. Namely, it is given by
the {\em fermionic formula} of~\cite{fjm}, and the generating function of
these characters for all $\alpha\in\Lambda_+$ is an eigenfunction of the
$q$-difference Toda integrable system. It would be interesting to find a
combinatorial relation between the monopole and fermionic formulas.

In the affine case, the zastava space $Z^\alpha_{\fg_\aff}$ was introduced 
in~\cite{bfg}. It is an irreducible affine algebraic variety containing a (dense
smooth symplectic) open subvariety $\oZ^\alpha_{\fg_\aff}$: the moduli space of
degree $\alpha$ based maps from the projective line $\BP^1$ to the Kashiwara
flag scheme ${\mathbf{Fl}}_{\fg_\aff}$. Contrary to the finite case, the open
subvariety $\oZ^\alpha_{\fg_\aff}$ is not affine, but only quasiaffine, and we
denote by $\ul\oZ{}^\alpha_{\fg_\aff}$ its affine closure. We do not know if the 
open embedding $\oZ^\alpha_{\fg_\aff}\hookrightarrow Z^\alpha_{\fg_\aff}$ extends to an
open embedding $\ul\oZ{}^\alpha_{\fg_\aff}\hookrightarrow Z^\alpha_{\fg_\aff}$:
it depends on the normality property of $Z^\alpha_{\fg_\aff}$ that is established
only for $\fg$ of types $A,C$~\cite{bf4,fr} at the moment (but is expected for 
all types). For an $A^{(1)},D^{(1)},E^{(1)}$ type quiver and an unframed quiver
gauge theory with $\alpha=\sum_{i\in\tilde{I}}(\dim V_i)\alpha_i$, we have an
isomorphism $\CM_C(\sG,\bN)\iso\ul\oZ{}^\alpha_{\fg_\aff}$~\cite{bfn3}. 
If $Z^\alpha_{\fg_\aff}$ is normal, this isomorphism extends to 
$\CM_C^+(\sG,\bN)\iso Z^\alpha_{\fg_\aff}$, and the fermionic formula for the
character $\chi(\BC[Z^\alpha])$ holds true.

Finally, for an arbitrary quiver $Q$ without edge loops we can consider an
unframed quiver gauge theory, and a coroot 
$\alpha:=\sum_{i\in Q_0}(\dim V_i)\alpha_i$ of the corresponding Kac-Moody Lie
algebra $\fg_Q$. The moduli space $\oZ^\alpha_{\fg_Q}$ of based maps from $\BP^1$
to the Kashiwara flag scheme ${\mathbf{Fl}}_{\fg_Q}$ was studied in~\cite{bfg}.
It is a smooth connected variety. We expect that it is quasiaffine, and its
affine closure $\ul\oZ{}^\alpha_{\fg_Q}$ is isomorphic to the Coulomb branch
$\CM_C(\sG,\bN)$. It would be interesting to find a modular interpretation of
$\CM_C^+(\sG,\bN)$ and its stratification into symplectic leaves. In the affine
case such an interpretation involves Uhlenbeck spaces $\CU_G^d(\BA^2)$.

\subsection{Multiplication and quantization}
\label{multiplication}
The multiplication of slices in the affine Grassmannian
$\ol\CW{}^\blambda_\bmu\times\ol\CW{}^{\blambda'}_{\bmu'}\to
\ol\CW{}^{\blambda+\blambda'}_{\bmu+\bmu'}$ was constructed in~\cite{bfn3} via 
multiplication of scattering matrices for singular monopoles (we learned
of its existence from T.~Dimofte, D.~Gaiotto and J.~Kamnitzer). The 
corresponding comultiplication $\BC[\ol\CW{}^{\blambda+\blambda'}_{\bmu+\bmu'}]\to
\BC[\ol\CW{}^\blambda_\bmu]\otimes\BC[\ol\CW{}^{\blambda'}_{\bmu'}]$ can not be seen
directly from the Coulomb branch construction of slices. However, its
quantization $\BC_\hbar[\ol\CW{}^{\blambda+\blambda'}_{\bmu+\bmu'}]\to
\BC_\hbar[\ol\CW{}^\blambda_\bmu]\otimes\BC_\hbar[\ol\CW{}^{\blambda'}_{\bmu'}]$
(recall from the end of~Section~\ref{general setup} that
$\BC_\hbar[\ol\CW{}^\blambda_\bmu]$ is the loop rotation equivariant Borel-Moore
homology of the corresponding variety of triples)
already can be realized in terms of Coulomb branches. The reason for this is
that the quantized Coulomb branch $\BC_\hbar[\ol\CW{}^\blambda_\bmu]$ is likely
to have a presentation by generators and relations of a {\em truncated shifted
Yangian} $Y^\blambda_\bmu\simeq\BC_\hbar[\ol\CW{}^\blambda_\bmu]$
of~\cite[Appendix~B]{bfn3}. Also, it seems likely that the comultiplication
of~\cite{fkprw} descends to a homomorphism $\Delta\colon 
Y^{\blambda+\blambda'}_{\bmu+\bmu'}\to Y^\blambda_\bmu\otimes Y^{\blambda'}_{\bmu'}$.
Finally, we expect that the desired comultiplication
$\BC[\ol\CW{}^{\blambda+\blambda'}_{\bmu+\bmu'}]\to
\BC[\ol\CW{}^\blambda_\bmu]\otimes\BC[\ol\CW{}^{\blambda'}_{\bmu'}]$ is obtained by
setting $\hbar=0$ in $\Delta$. 

Returning to the question of constructing the closed embeddings of slices
$\ol\CW{}_\bmu^{\blambda'}\hookrightarrow\ol\CW{}_\bmu^\blambda$ in terms of
Coulomb branches (see the beginning of~Section~\ref{warning}), we choose
dominant coweights $\bnu,\bnu',\bmu'$ such that $\bnu'+\bnu=\blambda',\
\bmu'+\bnu=\blambda'$, and set $\bmu'':=\bmu-\bmu'$. Then we have the
multiplication morphism $\ol\CW{}^\bnu_{\bmu''}\times\ol\CW{}^{\bnu'}_{\bmu'}\to
\ol\CW{}^\blambda_\bmu$, and we restrict it to $\ol\CW{}^\bnu_{\bmu''}=
\ol\CW{}^\bnu_{\bmu''}\times\{\bmu'\}\to\ol\CW{}^\blambda_\bmu$ where
$\bmu'\in\ol\CW{}^{\bnu'}_{\bmu'}$ is the fixed point. Then the desired closed
subvariety $\ol\CW{}_\bmu^{\blambda'}\subset\ol\CW{}_\bmu^\blambda$ is nothing but
the closure of the image of $\ol\CW{}^\bnu_{\bmu''}\to\ol\CW{}^\blambda_\bmu$.

The similar constructions are supposed to work for the slices in $\Gr_{G_\aff}$.
They are based on the comultiplication for {\em affine Yangians} constructed
in~\cite{gnw}.

\subsection{Affinization of KLR algebras}
We recall the setup of Section~\ref{quiver}, and set $W=0$ (no framing).
We choose a sequence $\bj=(j_1,\ldots,j_N)$ of vertices such that any vertex
$j\in Q_0$ enters $\dim V_j$ times; thus, $N=\dim V$. The set of all such 
sequences is denoted $\bJ(V)$. We choose a $Q_0$-graded flag
$V=V^0\supset V^1\supset\ldots\supset V^N=0$ such that $V^{n-1}/V^n$ is a
1-dimensional vector space supported at the vertex $j_n$ for any 
$n=1,\ldots,N$. It gives rise to the following flag of $Q_0$-graded lattices in
$V_\CK=V\otimes\CK\colon \ldots\supset L^{-1}\supset L^0\supset L^1\supset\ldots$,
where $L^{r+N}=zL^r$ for any $r\in\BZ;\ L^0=V_\CO$, and 
$L^n/L^N=V^n\subset V=L^0/L^N$ for any $n=1,\ldots,N$. Let $\sI_\bj\subset\sG_\CO$
be the stabilizer of the flag $L^\bullet$ (an Iwahori subgroup).
Then the $\sI_\bj$-module $\bN_\CO$ contains a submodule $\bN_\bj$ formed by  
the $\CK$-linear homomorphisms $b_e\colon V_{t(e),\CK}\to V_{h(e),\CK}$ such that
for any $e\in Q_1$ and $r\in\BZ,\ b_e$ takes $L^r_{t(e)}$ to $L^{r+1}_{h(e)}$. 

We consider the following version of the variety of triples:
$\CR_{\bj,\bj}:=\{[g,s]\in\sG_\CK\stackrel{\sI_\bj}{\times}\bN_\bj : gs\in\bN_\bj\}$,
cf.~\cite[Section~4]{bef,web}. Then the equivariant Borel-Moore homology
$\CH_{\bj,\bj}:=H_\bullet^{\BC^\times\ltimes\sI_\bj}(\CR_{\bj,\bj})$ forms an associative 
algebra with respect to a convolution operation. Moreover, if we take another 
sequence $\bj'\in\bJ(V)$ and consider 
$\CR_{\bj',\bj}:=\{[g,s]\in\sG_\CK\stackrel{\sI_\bj}{\times}\bN_\bj : 
gs\in\bN_{\bj'}\}$, then 
$\CH_{\bj',\bj}:=H_\bullet^{\BC^\times\ltimes\sI_{\bj'}}(\CR_{\bj,\bj'})$ forms a
$\CH_{\bj',\bj'}-\CH_{\bj,\bj}$-bimodule with respect to convolution, and we have
convolutions $\CH_{\bj'',\bj'}\otimes\CH_{\bj',\bj}\to\CH_{\bj'',\bj}$. In other words,
$\CH_V:=\bigoplus_{\bj,\bj'\in\bJ(V)}\CH_{\bj',\bj}$ forms a convolution algebra.

Furthermore, given $\bj_1\in\bJ(V),\ \bj_2\in\bJ(V')$, the concatenated
sequence $\bj_1\bj_2$ lies in $\bJ(V\oplus V')$, and one can define the 
morphisms $\CH_{\bj'_1,\bj_1}\otimes\CH_{\bj'_2,\bj_2}\to\CH_{\bj'_1\bj'_2,\bj_1\bj_2}$ 
summing up to a homomorphism $\CH_V\otimes\CH_{V'}\to\CH_{V\oplus V'}$. 

Similarly to the classical theory of Khovanov-Lauda-Rouquier algebras
(see a beautiful survey~\cite{rou} and references therein), we expect that
in case $Q$ has no loop edges, the categories of finitely generated graded
projective $\CH_V$-modules provide a categorification of the positive part 
$U_Q^{++}$ of the quantum toroidal algebra $U_Q$ (where $U_Q^{++}$ is defined as 
the subalgebra generated by the positive modes of the positive generators 
$e_{j,r}\colon j\in Q_0,\ r\geq0$).


\begin{thebibliography}{XXX}

\bibitem[BaKi]{baki} B.~Bakalov, A.~Kirillov, Jr., {\em Lectures on Tensor
categories and modular functors}, University Lecture Series {\bf 21}, American 
Mathematical Society, Providence, RI (2001), x+221pp.

\bibitem[Ba]{ba} V.~Baranovsky, {\em Uhlenbeck compactification as a functor},
Int.\ Math.\ Res.\ Not.\ IMRN (2015), no.~23, 12678--12712.

\bibitem[BD]{bd} A.~Beilinson, V.~Drinfeld, {\em Quantization of Hitchin's
integrable system and Hecke eigensheaves}, available at
http://www.math.uchicago.edu/$\tilde{\ }$mitya/langlands.html (2000).

\bibitem[Br]{br} T.~Braden, {\em Hyperbolic localization of intersection
cohomology}, Transform.\ Groups {\bf 8} (2003), 209--216.

\bibitem[BG]{bg} A.~Braverman, D.~Gaitsgory, {\em Crystals via the affine
Grassmannian}, Duke Math.\ J.\ {\bf 107} (2001), no.~3, 561--575.

\bibitem[BEF]{bef} A.~Braverman, P.~Etingof, M.~Finkelberg, {\em Cyclotomic
double affine Hecke algebras (with an appendix by Hiraku Nakajima and
Daisuke Yamakawa)}, arXiv:1611.10216.

\bibitem[BFGM]{bfgm} A.~Braferman, M.~Finkelberg, D.~Gaitsgory, I.~Mirkovi\'c,
{\em Intersection cohomology of Drinfeld's compactifications}, Selecta Math.\
(N.S.) {\bf 8} (2002) no.~3, 381--418; {\em Erratum}, Selecta Math.\
(N.S.) {\bf 10} (2004), 429--430.

\bibitem[BFG]{bfg} A.~Braverman, M.~Finkelberg, D.~Gaitsgory,
{\em Uhlenbeck spaces via affine Lie algebras}, Progr.\ Math.\ {\bf 244},
Birkh\"auser, Boston (2006), 17--135; {\em Erratum}, arxiv:0301176v4.

\bibitem[BF1]{bf1} A.~Braverman, M.~Finkelberg, {\em Pursuing the double affine
Grassmannian I: Transversal slices via instantons on $A_k$-singularities},
Duke Math.\ J.\ {\bf 152} (2010), no.~2, 175--206.

\bibitem[BF2]{bf2} A.~Braverman, M.~Finkelberg, {\em Pursuing the double affine
Grassmannian II: Convolution}, Advances in Math.\ {\bf 230} (2012), 414--432.

\bibitem[BF3]{bf3} A.~Braverman, M.~Finkelberg, {\em Pursuing the double affine
Grassmannian III: Convolution with affine zastava}, Moscow Math.\ J.\ {\bf 13}
(2013), no.~2, 233--265.

\bibitem[BF4]{bf4} A.~Braverman, M.~Finkelberg, {\em Semi-infinite Schubert 
varieties and quantum $K$-theory of flag manifolds}, J.\ Amer.\ Math.\ Soc.\ 
{\bf 27} (2014), no.~4, 1147--1168.

\bibitem[BFN1]{bfn1} A.~Braverman, M.~Finkelberg, H.~Nakajima, {\em Instanton
moduli spaces and $\CW$-algebras}, Ast\'erisque {\bf 385} (2016), vii+128pp.

\bibitem[BFN2]{bfn2} A.~Braverman, M.~Finkelberg, H.~Nakajima, {\em Towards a
mathematical definition of 3-dimensional $\CN=4$ gauge theories, II}, 
arXiv:1601.03586.

\bibitem[BFN3]{bfn3} A.~Braverman, M.~Finkelberg, H.~Nakajima, {\em Coulomb
branches of $3d\ \CN=4$ quiver gauge theories and slices in the affine
Grassmannian (with appendices by Alexander Braverman, Michael Finkelberg,
Joel Kamnitzer, Ryosuke Kodera, Hiraku Nakajima, Ben Webster, and Alex Weekes)},
arXiv:1604.03625.

\bibitem[BFN4]{bfn4} A.~Braverman, M.~Finkelberg, H.~Nakajima, {\em Ring 
objects in the equivariant derived Satake category arising from Coulomb 
branches (with appendix by Gus Lonergan)}, arXiv:1706:02112.

\bibitem[BK]{bk} A.~Braverman, D.~Kazhdan, {\em Representations of affine
Kac-Moody groups over local and global fields: a survey of some recent results},
European Congress of Mathematics, Eur.\ Math.\ Soc., Z\"urich (2013), 91--117.

\bibitem[BKP]{bkp} A.~Braverman, D.~Kazhdan, M.~Patnaik, {\em Iwahori-Hecke
algebras for $p$-adic loop groups}, Invent.\ Math.\ {\bf 204} (2016), no.~2,
347--442.

\bibitem[CHZ]{chz} S.~Cremonesi, A.~Hanany, A.~Zaffaroni, {\em Monopole
operators and Hilbert series of Coulomb branches of $3d\ \CN=4$ gauge theories},
JHEP {\bf 1401} (2014), 005.

\bibitem[CFHM]{cfhm} S.~Cremonesi, G.~Ferlito, A.~Hanany, N.~Mekareeya,
{\em Coulomb branch and the moduli space of instantons}, JHEP {\bf 1412} (2014),
103.

\bibitem[DG]{dg} V.~Drinfeld, D.~Gaitsgory, {\em On a theorem of Braden},
Transform.\ Groups {\bf 19} (2014), 313--358.

\bibitem[FJM]{fjm} B.~Feigin, E.~Feigin, M.~Jimbo, T.~Miwa, E.~Mukhin,
{\em Fermionic formulas for eigenfunctions of the difference Toda Hamiltonian}, 
Lett.\ Math.\ Phys.\ {\bf 88} (2009), no.~1-3, 39--77.

\bibitem[Fr]{Fr} I.~Frenkel, {\em Representations of affine Lie algebras,
Hecke modular forms and Korteweg de Vries type equations}, Lecture Notes in
Math.\ {\bf 933}, Springer, Berlin (1982), 71--110.

\bibitem[FK]{fk} M.~Finkelberg, D.~Kubrak, {\em Vanishing cycles on Poisson
varieties}, Funct.\ Anal.\ Appl.\ {\bf 49} (2015), no.~2, 135--141.

\bibitem[FiMi]{fm} M.~Finkelberg, I.~Mirkovi\'c, {\em Semi-infinite flags. I.
Case of global curve $\BP^1$}, Amer.\ Math.\ Soc.\ Transl.\ Ser.~2, {\bf 194},
Amer.\ Math.\ Soc., Providence, RI (1999), 81--112.

\bibitem[FeMa]{fema} E.~Feigin, I.~Makedonskyi, {\em Semi-infinite Pl\"ucker
relations and Weyl modules}, arXiv:1709.05674.

\bibitem[FKPRW]{fkprw} M.~Finkelberg, J.~Kamnitzer, K.~Pham, L.~Rybnikov,
A.~Weekes, {\em Comultiplication for shifted Yangians and quantum open Toda
lattice}, arXiv:1608.03331.

\bibitem[FR]{fr} M.~Finkelberg, L.~Rybnikov, {\em Quantization of Drinfeld
zastava in type $C$}, Algebr.\ Geom.\ {\bf 1} (2014), no.~2, 166--180.

\bibitem[GaRo]{garo} S.~Gaussent, G.~Rousseau, 
{\em Spherical Hecke algebras for Kac-Moody groups over local fields},
Ann.\ of Math.\ (2) {\bf 180} (2014), no.~3, 1051--1087. 

\bibitem[Gi]{gi} V.~Ginzburg, {\em Perverse sheaves on a loop group and
Langlands' duality}, arXiv:alg-geom/9511007.

\bibitem[GiRi]{gr} V.~Ginzburg, S.~Riche, {\em Differential operators on $G/U$
and the affine Grassmannian}, J.\ Inst.\ Math.\ Jussieu {\bf 14} (2015), no.~3,
493--575.

\bibitem[Gr]{gro} A.~Grothendieck, {\em Sur la classification des fibr\'es
holomorphes sur la sph\`ere de Riemann}, Amer.\ J.\ Math.\ {\bf 79} (1957), 
121--138.

\bibitem[GNW]{gnw} N.~Guay, H.~Nakajima, C.~Wendlandt, {\em Coproduct for
Yangians of affine Kac-Moody algebras}, arXiv:1701.05288.

\bibitem[Hi]{hi} T.~Hikita, {\em An Algebro-Geometric Realization of the
Cohomology Ring of Hilbert Scheme of Points in the Affine Plane}, 
Int.\ Math.\ Res.\ Not.\ IMRN (2017), no.~8, 2538--2561.

\bibitem[Ka1]{ka1} M.~Kashiwara, {\em The flag manifold of Kac-Moody Lie 
algebra}, Algebraic analysis, geometry and number theory (Baltimore, MD, 1988),
Johns Hopkins Univ.\ Press (1989), 161--190.

\bibitem[Ka2]{ka2} M.~Kashiwara, {\em On crystal bases}, CMS Conf.\ Proc.\
{\bf 16}, Amer.\ Math.\ Soc., Providence, RI (1995), 155--197.

\bibitem[KT]{kt} M.~Kashiwara, T.~Tanisaki, {\em Kazhdan-Lusztig conjecture
for affine Lie algebras with negative level}, Duke Math.\ J.\ {\bf 77} (1995),
21--62.

\bibitem[KWY]{kwy} J.~Kamnitzer, B.~Webster, A.~Weekes, O.~Yacobi, 
{\em Yangians and quantizations of slices in the affine Grassmannian}, 
Algebra and Number Theory {\bf 8} (2014), no.~4, 857--893.

\bibitem[KTWY]{ktwy} J.~Kamnitzer, P.~Tingley, B.~Webster, A.~Weekes, O.~Yacobi,
{\em Highest weights for truncated shifted Yangians and product monomial
crystals}, arXiv:1511.09131. 

\bibitem[Kr]{kr} V.~Krylov, {\em Integrable crystals and restriction to Levi
via generalized slices in the affine Grassmannian}, arXiv:1709.00391.

\bibitem[Lu]{lu} G.~Lusztig, {\em Singularities, character formulas, and a
$q$-analogue of weight multiplicities}, Ast\'erisque {\bf 101-102},
Soc.\ Math.\ France (1983), 208--229.

\bibitem[MS]{ms} G.~Moore, N.~Seiberg, {\em Classical and conformal field
theory}, Comm.\ Math.\ Phys.\ {\bf 123} (1989), 177--254.

\bibitem[MV]{mv} I.~Mirkovi\'c, K.~Vilonen, {\em Geometric Langlands duality 
and representations of algebraic groups over commutative rings}, Ann.\ of Math.\
(2) {\bf 166} (2007), no.~1, 95--143.

\bibitem[Na1]{na1} H.~Nakajima, {\em Geometric construction of representations
of affine algebras}, 
Proceedings of the International Congress of Mathematicians, {\bf I} 
(Beijing, 2002), Higher Ed.\ Press, Beijing (2002), 423--438.

\bibitem[Na2]{na2} H.~Nakajima, {\em Quiver varieties and branching}, 
SIGMA Symmetry Integrability Geom.\ Methods Appl.\ {\bf 5} (2009), paper 003.

\bibitem[Na3]{na3} H.~Nakajima, {\em Towards a mathematical definition of
Coulomb branches of 3-dimensional $\CN=4$ gauge theories, I}, Adv.\ Theor.\
Math.\ Phys.\ {\bf 20} (2016), no.~3, 595--669.

\bibitem[Na4]{na4} H.~Nakajima, {\em Questions on provisional Coulomb branches
of 3-dimensional $\CN=4$ gauge theories}, arXiv:1510.03908.

\bibitem[Na5]{na5} H.~Nakajima, {\em Introduction to a provisional
mathematical definition of Coulomb branches of 3-dimensional $\CN=4$
gauge theories}, arXiv:1706.05154.

\bibitem[NT]{nt} H.~Nakajima, Y.~Takayama, {\em Cherkis bow varieties and
Coulomb branches of quiver gauge theories of affine type $A$}, Selecta Math.\
(N.S.) {\bf 23} (2017), no.~4, 2553-2633.

\bibitem[R]{rou} R.~Rouquier, {\em Quiver Hecke algebras and 2-Lie algebras},
Algebra Colloq.\ {\bf 19} (2012), no.~2, 359--410.

\bibitem[Sl]{sl} W.~Slofstra, {\em A Brylinski filtration for affine Kac-Moody
algebras}, Adv.\ Math.\ {\bf 229} (2012), no.~2, 968--983.

\bibitem[SW]{sw} N.~Seiberg, E.~Witten, {\em Gauge dynamics and compactification
to three dimensions}, Adv.\ Ser.\ Math.\ Phys.\ {\bf 24}, World Sci.\ Publ.,
River Edge, NJ (1997), 333--366.

\bibitem[We]{web} B.~Webster, {\em Representation theory of the cyclotomic
Cherednik algebra via the Dunkl-Opdam subalgebra}, arXiv:1609.05494.

\bibitem[Wi]{wi} E.~Witten, {\em Geometric Langlands from six dimensions},
CRM Proc.\ Lecture Notes {\bf 50}, Amer.\ Math.\ Soc., Providence, RI (2010),
281--310.  
















\end{thebibliography}
\end{document}